\documentclass{amsart}

\usepackage{vmargin}		
\usepackage{fancyhdr}		
\usepackage[all]{xy}			
\usepackage{mathrsfs}		
\usepackage[T1]{fontenc}		
\usepackage{lmodern}
\usepackage{textcomp}
\usepackage[colorlinks=true,linkcolor=black,citecolor=black]{hyperref}				

\usepackage{amsmath}	
\usepackage{amssymb}	
\usepackage{amsfonts}	
\usepackage{mathrsfs}	
\usepackage{amsthm}	
\usepackage{graphicx}
\usepackage{url}
\usepackage{enumerate}

\theoremstyle{plain}  
\newtheorem{theorem}{Theorem}[subsection]
\newtheorem{lem}[theorem]{Lemma}
\newtheorem{prop}[theorem]{Proposition}
\newtheorem{hyp}[theorem]{Hypothesis}
\newtheorem{conj}[theorem]{Conjecture}
\newtheorem{cor}[theorem]{Corollary}

\newtheorem*{Conj}{Conjecture}

\newtheorem*{Prop}{Proposition}

\theoremstyle{definition}
\newtheorem{definition}[theorem]{Definition}

\newtheorem{remark}[theorem]{Remark}
\newtheorem{question}[theorem]{Question}

\makeatletter
\@addtoreset{theorem}{subsection}
\makeatother

\author{J\'er\^ome G\"artner}
\address{Insitut de Math\'ematiques de Jussieu\\
 Universit\'e Pierre et Marie Curie\\
 4 place Jussieu 75005 Paris France.} 
\email{jgartner@math.jussieu.fr}
\date{\today}

\begin{document}

\title[Darmon's points and quaternionic Shimura varieties]{Darmon's points and quaternionic Shimura varieties}

\begin{abstract}

In this paper, we generalize a conjecture due to Darmon and Logan (see \cite{DL} and \cite{Da}, chapter 8) in an adelic setting. We study the relation between our construction and Kudla's works on cycles on orthogonal Shimura varieties. This relation allows us to conjecture a Gross-Kohnen-Zagier theorem for Darmon's points.
\end{abstract}

\maketitle



\section{Introduction}

The theory of complex multiplication gives a collection of \textit{Heegner points} on elliptic curves over $\mathbf Q$, which are defined over class fields of imaginary quadratic fields. These points allowed to prove Birch and Swinnerton-Dyer's conjecture over~$\mathbf Q$ for analytic rank 1 curves, thanks to the work of Gross-Zagier and Kolyvagin.

Let us briefly recall the construction of Heegner points. If $E$ is an elliptic curve over $\mathbf Q$ then we know that $E$ is modular. Let $N$ be the conductor of $E$. There exists a modular form $f\in S_2(N)$ such that $L(E,s)=L(f,s)$. Denote by $\Phi_N:\Gamma_0(N)\backslash \mathcal H\longrightarrow E(\mathbf C)$ the modular uniformization which is obtained by taking the composition of the map ${z_0\in\mathcal H\mapsto c\int_{i\infty}^{z_0} 2\pi i f(z)\mathrm{d}z}$ (here $c$ denotes the Manin constant) with the Weierstrass uniformization. Let $z_0\in\mathcal H\cap K$, where $K/\mathbf Q$ is an imaginary quadratic field. A Heegner point is given essentially by $2\pi i\int_{i\infty}^{z_0} f(z)\mathrm{d}z$ modulo periods of $f$. It is the Abel-Jacobi image of $z_0$ in $\mathbf C/\Lambda_E\simeq E(\mathbf C)$. The theory of complex multiplication shows that these points are defined over class fields of $K$.

In \cite{Da}, Darmon gives a conjectural construction of \textit{Stark-Heegner points}, which is a generalization of classical Heegner points. These points should help us to understand, on one hand the Birch and Swinnerton-Dyer conjecture, on the other hand Hilbert's twelfth problem.

In more concrete terms, assume that $F$ is a totally real number field of narrow class number 1. Let $\tau_j$ be its archimedean places, and $K/F$  some quadratic ``ATR'' extension (\text{i.e.} $K$ has exactly one complex place). Darmon defines a collection of points on elliptic curves $E/F$ which are expected to be defined over class fields of $K$. In this case, the (conjectural, but partially proved by Skinner - Wiles) modularity of $E$ gives the existence of a Hilbert modular form $f$ on $\mathcal H^r$ whose periods appear as a tensor product of periods of $E_{\tau_j}=E\otimes_{F,\tau_j}\mathbf C$. The construction explained in \cite{DL} can be seen as an exotic Abel-Jacobi map.

In this paper, we  generalize Darmon's contruction by removing the hypothesis ``ATR'' on $K$ (but we assume that $K$ is not CM) and the technical hypothesis that $F$ has narrow class number~1. We replace the Hilbert modular variety used in the ``ATR'' case by a general quaternionic Shimura variety and define a suitable Abel-Jacobi map. We are able to specify the invariants of the quaternion algebra using local epsilon factors and to give a conjectural Gross-Zagier formula for these points. We conclude the paper by establishing a relation to Kudla's study of cycles on orthogonal Shimura varieties, in order to give a Gross-Kohnen-Zagier type conjecture.

Let us summarize the main construction of this paper. Let $F$ be a totally real field of degree $d$ and let $\tau_1,\dots,\tau_d$ be its archimedean places. Fix $r\in\{2,\dots, d\}$, and and a quadratic extension $K/F$ such that the set of archimedean places of $F$ that split completely in $K$ is $\{\tau_2,\dots,\tau_r\}$. Let $B/F$ be a quaternion algebra which splits at $\tau_1,\dots,\tau_r$ and ramifies at $\tau_{r+1},\dots,\tau_d$. Let $G=\mathrm{Res}_{F/\mathbf Q} B^{\times}$. We will denote by $\mathrm{Sh}_H(G,X)$ the quaternionic Shimura variety of level $H$ (a compact open subgroup of $G(\mathbf A_f)$) whose complex points are given by
$$\mathrm{Sh}_H(G,X)(\mathbf C)=G(\mathbf Q)\backslash (\mathbf C\smallsetminus\mathbf R)^r\times G(\mathbf A_f)/H .$$

Fix an embedding $q:K\hookrightarrow B$. There is an action of $(K\otimes \mathbf R)^{\times}_+/(F\otimes\mathbf R)^{\times}$ on $(\mathbf C\smallsetminus\mathbf R)^r$. By considering a suitable orbit of this action, we obtain a real cycle $T_b$ of dimension $r-1$ on $\mathrm{Sh}_H(G,X)(\mathbf C)$. Using the theorem of Matsushima and Shimura, we deduce that there exists an $r$-cycle $\Delta_b$ on $\text{Sh}_H(G,X)(\mathbf C)$ such that $\partial \Delta_b$ is an integral multiple of $\mathscr T_b$.  

Let $E/F$ be an elliptic curve, assumed modular, {i.e.}, there exists a Hilbert modular eigenform $\tilde{\varphi}$ satisfying $L(E,s)=L(\tilde{\varphi},s)$. We will assume that this form corresponds to an automorphic form $\varphi$ on $B$ by the Jacquet-Langlands correspondence. There exists a holomorphic differential form $\omega_{\varphi}$ of degree $r$ on $\mathrm{Sh}_H(G,X)(\mathbf C)$ naturally attached to $\varphi$. In general, the set of periods of $\omega_{\varphi}$ is a dense subset of $\mathbf C$. Fix some character $\beta$ of the set of connected components of $(K\otimes\mathbf R)_+^{\times}/(F\otimes\mathbf R)^{\times}$. Following Darmon we define a modified differential
form $\omega_{\varphi}^{\beta}$ whose periods are, assuming Yoshida's period conjecture, a lattice, homothetic to some sublattice of the Neron lattice of $E$. 

The image of (a suitable
multiple of) the complex number $\int_{\Delta_b} \omega_{\varphi}^{\beta}$  in $\mathbf C/\Lambda_E$ is independent of the
choice of $\Delta_b$. Hence it defines by Weierstrass uniformization a point $P_b^{\beta}$ in $E(\mathbf C)$. We conjecture

\begin{Conj}[\ref{conjectureprincipale}]
 
$P_b^{\beta}=\Phi\left(\int_{\Delta_b}\omega_{\varphi}^{\beta}\right)\in E(\mathbf C)$ lies in $E(K^{ab})$ and 
 $$\forall a\in\mathbf A_{K}^{\times}\qquad\mathrm{rec}_K(a)P_b^{\beta}=\beta(a_{\infty})P_{q_{\mathbf A}(a)b}^{\beta}.$$

\end{Conj}

Let us assume this conjecture is true and denote by $K_b^+$ the field of definition of $P_b^{\beta}$. Let $\pi=\pi(\varphi)$ be the automorphic representation generated by $\varphi$; fix a character $\chi:\mathrm{Gal}(K_b^+/K)\rightarrow\mathbf C^{\times}$. Denote by $\varepsilon(\pi\times\chi,\frac 12)$ the sign in the functional equation of the Rankin-Selberg $L$-function $L(\pi\times\chi,s)$ and by $\eta_K:F_{\mathbf A}^{\times}/F^{\times}\mathrm{N}_{K/F}(K_{\mathbf A}^{\times})\rightarrow\{\pm 1\}$ the quadratic character of $K/F$. The following proposition proves that $B$ is uniquely determined by $K$ and the isogeny class of $E/F$.

\begin{Prop}[\ref{determinationdeB}]
 Let $b\in \widehat{B}^{\times}$ and assume conjecture \ref{conjectureprincipale}. If
$$e_{\overline{\chi}}(P_b^{\beta})=\sum_{\sigma\in\mathrm{Gal}(K_b^+/K)} \chi(\sigma)\otimes P_b^{\beta}\in E(K_b^+)\otimes \mathbf Z[\chi]$$
\noindent is not torsion, then :
$$\forall v\nmid\infty\qquad \eta_{K,v}(-1)\varepsilon(\pi_v\times\chi_v,\frac 12)=\mathrm{inv}_v(B_v)\quad\mathrm{and}\quad \varepsilon(\pi\times\chi,\frac 12)=-1.$$ 
\end{Prop}

The last part of this paper is focused on a conjecture in the spirit of the Gross-Kohnen-Zagier theorem. Assume that $E(F)$ has rank 1. Denote by $P_0$ some generator modulo torsion. For each totally positive
$t\in O_F$ such that $(t)$ is square free and prime to $d_{K/F}$, denote by $K[t]$ the quadratic extension $K[t]=F(\sqrt{-D_0t})$, where $D_0\in F$ satisfies $\tau_j(D_0)>0$ if and only if $j\in\{1,r+1,\dots,d\}$.
Let $P_{t,1}$ be Darmon's point obtained for $K[t]$ and $b=1$, and set $$P_t=\mathrm{Tr}_{K[t]_1^+/F} P_{t,1}.$$
\noindent The point $P_t$ is in $E(F)$ and there exists some integer $[P_t]\in\mathbf Z$ such that $P_t=[P_t]P_0.$
In the spirit of conjecture 5.3 of \cite{DT} we conjecture that :

\begin{Conj}[\ref{lastconj}]
There exists a Hilbert modular form $g$ of level 3/2 such that the $[P_t]$s are proportional to some Fourier coefficients of $g$.
\end{Conj}

In our attempt to adapt Yuan, Zhang and Zhang's proof in the CM case \cite{YZZ} to prove this conjecture, we obtained a relation between Darmon's points and Kudla's program, see Proposition~5.5.3.2.

\proof[Acknowledgments] This work grew out of the author's thesis at University Paris 6. The author is grateful to J. Nekov\'a\v r for his constant support during this work, and to C. Cornut for many useful conversations.

\section{Quaternionic Shimura varieties}
\addtocounter{subsection}{1}

In this section we recall some properties of Shimura varieties associated to quaternion algebras. The standard references are Reimann's book \cite{R} and \cite{Mi2}. The content of this section is more or less the transcription to Shimura varieties of what is done for curves in \cite{CV} and \cite{N}.

Let $F$ be a totally real field of degree $d=[F:\mathbf{Q}]$ and $\tau_1,\dots,\tau_d$ its archimedean places. Denote by $\overline{\mathbf{Q}}\subset\mathbf{C}$ the algebraic closure of $\mathbf{Q}$ in $\mathbf{C}$ so 
$\tau_j:F\lhook\joinrel\relbar\joinrel\rightarrow\overline{\mathbf{Q}}.$ Fix $r\in\{2,\dots,d\}$ and a finite set $S_B$ of non-archimedean primes satisfying $$|S_B|\equiv d-r\ \mathrm{mod}\ 2.$$ 
\noindent Let $B$ be the unique quaternion algebra over $F$ ramified at the set
$$\mathrm{Ram}(B)=\{\tau_{r+1},\dots,\tau_d\}\cup S_B.$$

\noindent For each $j\in\{1,\dots,d\}$ we put $B_{\tau_j}=B\otimes_{F,\tau_j}\mathbf{R}$. It is not necessary but more convenient to fix for each $j\in\{\tau_1,\dots,\tau_r\}$ an $\mathbf{R}$-algebra isomorphism $$B_{\tau_j}\overset{\sim}{\longrightarrow} M_2(\mathbf{R}).$$
\noindent The constructions given in this paper are independent on the choice of these isomorphisms, as in the author's PhD thesis \cite{JG}.

Let $G$ be the algebraic group over $\mathbf{Q}$ satisfying $G(A)=(B\otimes_{\mathbf{Q}}A)^{\times}$ for every commutative $\mathbf{Q}$-algebra $A$. We will denote by $\mathrm{nr}:G(A)\longrightarrow(F\otimes_{\mathbf{Q}}A)^{\times}$ the reduced norm and by $Z$ the center of $G$. For $j\in\{1,\dots,d\}$ let $G_j$ be the algebraic group over $\mathbf{R}$ given by $G_j=G\otimes_{F,\tau_j}\mathbf{R}$; thus $G_{\mathbf R}$ decomposes as $G_1\times\dots\times G_d$. For any abelian group $A$, denote by $\widehat{A}$ the group $A\otimes\widehat{\mathbf{Z}}$.

Let $X$ be the $G(\mathbf{R})$-conjugacy class of the morphism 
$h:\mathbf{S}=\mathrm{Res}_{\mathbf{C}/\mathbf{R}}(\mathbf{G}_{m,\mathbf{C}})\longrightarrow
G(\mathbf{R})=G_1(\mathbf{R})\times\dots\times G_d(\mathbf{R})$ defined by
$$x+iy\longmapsto\left(\underbrace{\left(\begin{matrix} x&y\\-y&x\end{matrix}\right),\dots,
\left(\begin{matrix} x&y\\-y&x\end{matrix}\right)}_{r\ \mathrm{times}},\underbrace{1,\dots,1}_{d-r\ \mathrm{times}}\right).$$

\noindent The set $X$ has a natural complex structure \cite{Mi1} and the following map is an holomorphic isomorphism between $X$ and $(\mathbf{C}\smallsetminus\mathbf{R})^
r$ :
$$ghg^{-1}\longmapsto g\cdot (i,\dots,i)=\left(\frac{a_1 i +b_1}{c_1 i +d_1},\dots,\frac{a_r i +b_r}{c_ri+d_r}\right),$$
\noindent where $g=(g_1,\dots,g_d)\in G(\mathbf{R})$ and for $j\in\{1,\dots,r\}$ $g_j$ is identified with $\begin{pmatrix} a_j&b_j\\ c_j&d_j\end{pmatrix}$.

\paragraph{Quaternionic Shimura varieties}  

Let $H$ be an open-compact subgroup of $\widehat{B}^{\times}$. The quaternionic Shimura varieties considered in this paper are algebraic varieties $\mathrm{Sh}_H(G,X)$ whose complex points are given by
$$\mathrm{Sh}_H(G,X)(\mathbf{C})=B^{\times}\backslash(X\times\widehat{B}^{\times}/H)),$$
\noindent where the left-action of $B^{\times}$ and the right-action of $H$ are given by
$$\forall k\in B^{\times}\ \forall h\in H\ \forall (x,b)\in X\times \widehat{B}^{\times}\qquad k\cdot(x,b)\cdot h=(kx,kbh).$$

\noindent Such Shimura varieties are defined over some number field called the reflex field. In our case this number field is
$$F'=\mathbf{Q}\left(\sum_{j=1}^r\tau_j(\alpha),\ \alpha\in F\right)\subset\overline{\mathbf{Q}}\subset\mathbf{C}.$$

\noindent We will denote by $[x,b]_H$ the element of $\mathrm{Sh}_H(G,X)(\mathbf{C})$ represented by $(x,b)$ and by 
$[x,b]_{H\widehat{F}^{\times}}$ the corresponding element of the modified variety 
$\mathrm{Sh}_H(G/Z,X)(\mathbf{C})=B^{\times}\backslash(X\times\widehat{B}^{\times}/HZ))$.

\begin{remark}
 The complex Shimura varieties are compact whenever $B\neq M_2(F)$. The Hilbert modular varieties used by Darmon in \cite{Da} chapter 7 and 8 are obtained when $B=M_2(F)$ and $r=d$.
\end{remark}

The Shimura varieties form a projective system $\{\mathrm{Sh}_H(G,X)\}_H$ indexed by open compact subgroups in $\widehat{B}^{\times}$. The transition maps $\mathrm{pr}:\mathrm{Sh}_H(G,X)\rightarrow \mathrm{Sh}_{H'}(G,X)$ are defined on complex points by
$$[x,b]_H\rightarrow[x,b]_{H'}.$$

\noindent There is an action of $\widehat{B}^{\times}$ on the projective system $\{\mathrm{Sh}_H(G,X)\}_H$. The right multiplication by $g\in\widehat{B}^{\times}$ induces an isomorphism  $[\cdot g]: \{\mathrm{Sh}_H(G,X)\}_H\overset{\sim}{\longrightarrow}\{\mathrm{Sh}_H(G,X)\}_{g^{-1}Hg}$, defined on complex points by
$$[\cdot g][x,b]_H=[x,bg]_{g^{-1}Hg}.$$

\paragraph{Complex conjugation} 

Fix $j\in\{1,\dots,r\}$.
Let $h_j:\mathbf{S}\rightarrow G_{j,\mathbf{R}}$ be the morphism obtained by composing $h$ with the $j$-th projection $G_{\mathbf{R}}\rightarrow G_{j, \mathbf{R}}$ and $X_j$ the $G_j(\mathbf{R})$-conjugacy class of $h_j$. For $x_j=g_jh_jg_j^{-1}\in X_j$, the set $\mathrm{Im}(g_jh_jg_j^{-1})$ is a maximal anisotropic $\mathbf{R}$-torus in $G_{j,\mathbf{R}}$. The map $\ell_j : x_j\mapsto \mathrm{Im}(x_j)$ satisfies $|\ell_j^{-1}(\ell_j(x_j))|=2$, thus there exists a unique antiholomorphic and $G_{j,\mathbf{R}}$-equivariant involution $$t_j:X_j\longrightarrow X_j$$

\noindent such that $$\forall x_j\in X_j\qquad \ell_j^{-1}(\ell_j(x_j))=\{x_j,t_j(x_j)\}.$$

\noindent More precisely, under the identification $X_j\overset{\sim}{\longrightarrow}\mathbf{C}\smallsetminus\mathbf{R}$, the map $\ell_j$ satisfies $\ell_j(x+iy)=\left\{\begin{pmatrix}
x&y\\-y&x\end{pmatrix}\right\}$ and $\ell_j^{-1}(\ell_j(x+iy))=\{x+iy,x-iy\}$. Note that the map $t_j$ can be extended to complex points of the Shimura varieties by $t_j([x,b]_H)=[t_j(x),b]_H$; $t_j$ acts trivially on $X_k$ for $k\neq j$.

\paragraph{Differential forms} In this section we recall some facts concerning differential forms on Shimura varieties. We will denote by $\Omega_H=\Omega_{H/F'}$ the sheaf of differentials of degree $r$ on $\mathrm{Sh}_H(G,X)$ and by $\Omega_H^{\mathrm{an}}$ the sheaf of holomorphic $r$-differentials on $\mathrm{Sh}_H(G,X)(\mathbf{C})$,  provided that $\mathrm{Sh}_H(G,X)$ is smooth. Recall that the GAGA principle gives us the following isomorphism between global sections
$$\Gamma(\mathrm{Sh}_H(G,X) ,\Omega_H)\otimes_{F'}\mathbf{C}\overset{\sim}{\longrightarrow}\Gamma(\mathrm{Sh}_H(G,X)(\mathbf{C}),\Omega_H^{\mathrm{an}}).$$

\noindent Notice that in general, $\mathrm{Sh}_H(G,X)$ is not smooth. In this last case we will fix some integer $n\geq 3$ such that for each $\mathfrak p$ in $\mathrm{Ram}(B)$ we have $\mathfrak p \nmid n\mathcal O_F$ and for each $v\mid n\mathcal O_F$ isomorphisms $\iota_v:B_v\overset{\sim}{\rightarrow}M_2(F_v)$. The group
$$H'=\left\{(h_v)\in H,\ \mathrm{s.t.}\ \forall v\mid n\mathcal O_F \ h_v\equiv\begin{pmatrix}
1&0\\0&1\end{pmatrix}\ \mathrm{mod}\ n\mathcal O_{F_v}\right\}$$
\noindent is of finite index in $H$ and $\mathrm{Sh}_{H'}(G,X)$ is smooth. The map $\mathrm{Sh}_{H'}(G,X)\rightarrow \mathrm{Sh}_H(G,X)$ is a finite covering. We define $\Omega_H=\frac{1}{[H:H']}\sum_{\sigma\in H/H'}\sigma\Omega_{H'}=(\Omega_{H'})^H$. By abuse of language, we shall call an element of $\Gamma(\Omega_H)=\Gamma(\mathrm{Sh}_H(G,X),\Omega_H)=(\sum_{\sigma\in H/H'}\sigma)\Gamma(\mathrm{Sh}_{H'}(G,X),\Omega_{H'})$ a global $r$-form on $\mathrm{Sh}_H(G,X)$. Remark that the space of global holomorphic $r$-forms $\varinjlim_H\Gamma(\Omega_H^{\mathrm{an}})$ is equipped with a canonical action of $\widehat{B}^{\times}$ given by pull-backs $[\cdot g]^*$.

Let $\varepsilon\in\{\pm 1\}^r$ and denote by $\Gamma((\Omega_H^{\mathrm an})^{\varepsilon})$ the space of $r$-forms on $\mathrm{Sh}_H(G,X)(\mathbf{C})$ which are holomorphic (resp. anti-holomorphic) in $z_j$ if $\varepsilon_j=+1$ (resp. if $\varepsilon_j=-1$). The maps $t_j$ pulled-back on $\Gamma((\Omega_H^{\mathrm an})^{\varepsilon})$ satisfy
$$t_j^*:\Gamma((\Omega_H^{\mathrm an})^{\varepsilon})\longrightarrow\Gamma((\Omega_H^{\mathrm an})^{\varepsilon'})$$

\noindent where $\varepsilon'_k=\varepsilon_k$ for $k\neq j$ and $\varepsilon'_j=-\varepsilon_j$. 

\noindent When $\sigma\in\prod_{j=2}^r\{\pm 1\}$ we will define $e_j\in\{0,1\}$ by $\sigma_j=(-1)^{e_j}$ and $t_{\sigma}^*$ by $\prod_{j=2}^r(t_j^*)^{e_j}$.
Let $\beta:\prod_{j=2}^r\{\pm 1\}\rightarrow\{\pm 1\}$ be a character and $\omega\in\Gamma(\Omega_H^{\mathrm{an}})$. We shall denote by $\omega^{\beta}$ the element $\omega^{\beta}=\sum_{\sigma\in\{\pm 1\}^{r-1}}\beta(\sigma)t_{\sigma}^*(\omega)$ of $\bigoplus_{\varepsilon}\Gamma((\Omega_H^{\mathrm{an}})^{\varepsilon})$.

\paragraph{Automorphic forms} Let $S_2^H$ be the space $S_{2,\dots,2,0,\dots,0}^H(B_{\mathbf{A}}^{\times})$ of functions 
$\varphi:B_{\mathbf{A}}^{\times}\simeq G(\mathbf{R})\times\widehat{B}^{\times}\longrightarrow\mathbf{C}$
satisfying the following properties :
\begin{enumerate}
 \item $\forall g\in B^{\times}\ \forall b\in B_{\mathbf{A}}^{\times}\qquad \varphi(gb)=\varphi(b)$,
\item $\forall g\in(\mathbf{R}^{\times})^r\times G_{r+1}(\mathbf{R})\times\dots\times G_d(\mathbf{R})\subset G(\mathbf{R})\ \forall b\in B_{\mathbf{A}}^{\times}\qquad \varphi(bg)=\varphi(b)$,
\item $\forall h\in H\ \forall b\in B_{\mathbf{A}}^{\times}\qquad \varphi(bh)=\varphi(b)$,
\item $\forall g\in B_{\mathbf{A}}^{\times}\ \forall (\theta_1,\dots,\theta_r)\in\mathbf{R}^r $ 
\begin{multline*}\varphi\left( g\left[\left(\begin{array}{cc}\cos\theta_1 & -\sin\theta_1 \\\sin\theta_1 & \cos\theta_1\end{array}\right),\dots
,\left(\begin{array}{cc}\cos\theta_r & -\sin\theta_r \\\sin\theta_r & \cos\theta_r\end{array}\right),1,\dots,1\right]\right)=\\
e^{-2i\theta_1}\times\dots\times e^{-2i\theta_r} \varphi(g),
\end{multline*}
\item For all $g\in B_{\mathbf{A}}^{\times}$, the map
$$(x_1+iy_1,\dots,x_r+iy_r)\mapsto\frac 1{y_1\dots y_r} \varphi\left(g\left[\left(\begin{array}{cc}y_1 & x_1 \\0 & 1\end{array}\right),\dots,\left(\begin{array}{cc}y_r & x_r \\0 & 1\end{array}\right),1,\dots,1\right]\right)$$
\noindent is holomorphic on $\mathcal{H}^r$ where $\mathcal H$ denotes the  Poincar\'e upper-half plane.
\end{enumerate}

\noindent Remark that we do not need any assumption to obtain cuspidal forms as $B$ will be assumed to differ from $M_2(F)$.

\noindent There is an action of $\widehat{B}^{\times}$ on $S_2=\bigcup_H S_2^H$ defined by
$$\forall g\in\widehat{B}^{\times},\ \forall\varphi\in S_2,\ \forall x\in B_{\mathbf{A}}^{\times}\qquad g\cdot\varphi(x)=\varphi(xg);$$
\noindent thus $S_2^H$ is the space of $H$-invariant functions in $S_2$.

By modifying the properties 4 and 5 above we obtain the following new definition :

\begin{definition}
 Let $\varepsilon:\{\tau_1,\dots,\tau_r\}\longrightarrow\{\pm 1\}$ and $\varepsilon_i=\varepsilon(\tau_i)$. The space $(S_2^{\varepsilon})^H$ is the space of maps $\varphi:B_{\mathbf{A}}^{\times}\simeq G(\mathbf{R})\times\widehat{B}^{\times}\longrightarrow\mathbf{C}$
\noindent satisfying 1-3 above and
\begin{enumerate}
 \item[4'.]  for all $g\in B_{\mathbf{A}}^{\times}$ and $(\theta_1\dots \theta_r)\in\mathbf{R}^r$  
\begin{align*}
\varphi\left( g\left(\left(\begin{array}{cc}\cos\theta_1 & -\sin\theta_1 \\\sin\theta_1 & \cos\theta_1\end{array}\right),\dots,\right.\right.&\left.\left.\left(\begin{array}{cc}\cos\theta_r & -\sin\theta_r \\\sin\theta_r & \cos\theta_r\end{array}\right),1,\dots,1\right)\right)\\
&=e^{-2i\varepsilon_1\theta_1}\times\dots\times e^{-2i\varepsilon_r\theta_r} \varphi(g)
\end{align*}
\item[5'.] for all $g\in B_{\mathbf{A}}^{\times}$ the map 
$$(x_1+iy_1,\dots,x_r+iy_r)\mapsto\frac 1{y_1\dots y_r} \varphi\left(g\left(\left(\begin{array}{cc}y_1 & x_1 \\0 & 1\end{array}\right),\dots,\left(\begin{array}{cc}y_r & x_r \\0 & 1\end{array}\right),1,\dots,1\right)\right)$$
\noindent is holomorphic (resp. anti-holomorphic) in $z_j=x_j+iy_j\in\mathcal H$ if $\varepsilon_j=1$ (resp. $\varepsilon_j=-1$).
\end{enumerate}

\end{definition}

\noindent We will denote by $S_2^{\hat{F}^{\times}}$ (resp. $(S_2^{\varepsilon})^{\hat{F}^{\times}}$) the space of elements in $S_2$ (resp. $S_2^{\varepsilon}$) which are $\hat{F}^{\times}$-invariant.

We are now able to affirm the existence of relations between $r$-forms on $\mathrm{Sh}_H(G,X)(\mathbf{C})$ and automorphic forms :
\begin{prop}
There exist bijections compatible with the $\widehat{B}^{\times}$-action between the following spaces :
\begin{center}
\begin{tabular}{lll}
$\Gamma(\Omega_H^{\mathrm{an}})$ &and & $S_2^H$\\
$\Gamma((\Omega_H^{\mathrm{an}})^{\varepsilon})$ &and & $(S_2^{\varepsilon})^H$\\
$\Gamma(\mathrm{Sh}_H(G/Z,X)(\mathbf{C}),(\Omega_H^{\mathrm{an}})^{\varepsilon})$ &and & $(S_2^{\varepsilon})^{H\widehat{F}^{\times}}$
\end{tabular}
\end{center}

\end{prop}

This statement is completely analogous to section 3.6 of \cite{CV}, see \cite{JG}, Propositions 1.2.2.4 and 1.2.2.5 for more details.

\paragraph{Matsushima-Shimura theorem}

The decomposition of the cohomology of quaternionic Shimura varieties  given by Matsushima-Shimura theorem will be usefull in the following sections. Let us recall this result when $B\neq M_2(F)$  \cite{MS} and \cite{F}. Denote by $h_F^+$ the narrow class number of $F$.

\begin{theorem}
\label{theoremedematsushimashimura}
 Let $m\in\{0,\dots,2r\}$. We have the following decomposition :
$$H^m(\mathrm{Sh}_H(G,X)(\mathbf{C}),\mathbf{C})\simeq\left\{\begin{array}{ll}
                                 \left(\mathrm{Vect} \bigwedge_{\substack{i\in a\subset\{1,\dots,r-1\}\\ |a|=m/2}} \frac{\mathrm{d} z_i\wedge\mathrm{d} \overline{z_i}}{y_i^2}\right)^{s}&\mathrm{if}\ m\neq r\\[4mm]
\left(\mathrm{Vect} \bigwedge_{\substack{i\in a\subset\{1,\dots,r-1\}\\ |a|=m/2}} \frac{\mathrm{d} z_i\wedge\mathrm{d} \overline{z_i}}{y_i^2}\right)^{s}\oplus\bigoplus_{\varepsilon\in\{\pm 1\}^r}(S_2^{\varepsilon})^H &\mathrm{if}\ m=r
                                \end{array}\right.
$$
\noindent and
$$H^m(\mathrm{Sh}_H(G/Z,X)(\mathbf{C}),\mathbf{C})\simeq\left\{\begin{array}{ll}
                                 \left(\mathrm{Vect} \bigwedge_{\substack{i\in a\subset\{1,\dots,r-1\}\\ |a|=m/2}} \frac{\mathrm{d} z_i\wedge\mathrm{d} \overline{z_i}}{y_i^2}\right)^{s'}&\mathrm{if}\ m\neq r\\[4mm]
\left(\mathrm{Vect} \bigwedge_{\substack{i\in a\subset\{1,\dots,r-1\}\\ |a|=m/2}} \frac{\mathrm{d} z_i\wedge\mathrm{d} \overline{z_i}}{y_i^2}\right)^{s'}\oplus\bigoplus_{\varepsilon\in\{\pm 1\}^r}(S_2^{\varepsilon})^{H\hat{F}^{\times}} &\mathrm{if}\ m=r,
                                \end{array}\right.
$$
\noindent where $s$ (resp. $s'$) is the number of connected components of $\mathrm{Sh}_H(G,X)(\mathbf C)$ (resp. of $\mathrm{Sh}_H(G/Z,X)(\mathbf C)$).
\end{theorem}





\section{Periods}

\subsection{Yoshida's conjecture}

Let $E/F$ be an elliptic curve, assumed modular in the sense that there exists
a cuspidal, parallel weight two Hilbert modular form $\tilde{\varphi}\in S_2(\mathrm{GL}_2(F_{\mathbf{A}}))$ satisfying $L(E,s)=L(\tilde{\varphi},s)$. We shall assume that the automorphic representation generated by $\tilde{\varphi}$ is obtained by the Jacquet-Langlands correspondence from
$\varphi\in S_2^{H\widehat{F}^{\times}}(B_{\mathbf{A}}^{\times})$.

Denote by $\pi=\pi_{\infty}\otimes\pi_f$ the automorphic representation of $B_{\mathbf{A}}^{\times}/F_{\mathbf{A}}^{\times}$ generated by $\varphi$. We shall assume until section \ref{sectionmultiplicite}, only for simplicity, that $\mathrm{dim}\ \pi_f^H=1$. 

Let $M=h^1(E)$ be the motive over $F$ with coefficients in $\mathbf{Q}$ associated to $E$. Yoshida \cite{Y} conjectures the existence of a rank $2^r$ motive $M'$ over the reflex field $F'$, with coefficients in $\mathbf{Q}$, satisfying $M'=\bigotimes_{\{\tau_1,\dots,\tau_r\}}\mathrm{Res}_{F/F'} M$. This motivic conjecture is the following :

\begin{conj}[Yoshida, \cite{Y}]

The motive $M'$ over $F'$ is isomorphic to the motive associated to the part $H^*(\mathrm{Sh}_{H\widehat{F}^{\times}}(G,X))^{(E)}$ of the cohomology for which Hecke eigenvalues are the same as $E$.

\end{conj}

While looking at the $\ell$-adic realization, this conjecture is in fact the Langlands cohomological conjecture. This case is known, up to semi-simplification, thanks to Brylinski and Labesse in the case $B=\mathrm{M}_2(F)$ \cite{BL}, Langlands in the case $B\neq \mathrm{M}_2(F)$ for primes of good reduction, \cite{L} and Reimann (- Zink) \cite{R, RZ} for a more general cases.

Recall the following decompositions given by Yoshida in \cite{Y} section 5.1, when we focus on $\tau':F'\hookrightarrow\mathbf{C}$ induced by $\tilde{\tau'}:\overline{\mathbf{Q}}\hookrightarrow\mathbf{C}$.

\paragraph{Betti cohomology} There exists an isomorphism of $\mathbf{Q}$-vector spaces
$$\mathscr I : M_{\mathrm B}'\overset{\sim}{\longrightarrow}\bigotimes_{j=1}^r M_{B,\tau_j}$$

\paragraph{de Rham cohomology} The map
$$\mathscr J : M_{\mathrm{dR}}'\overset{\sim}{\longrightarrow}\left(\bigotimes_{j=1}^r\left(M_{\mathrm{dR}}\otimes_{F,\tau_j} \overline{\mathbf{Q}}\right)\right)^{\mathrm{Gal}(\overline{\mathbf{Q}}/F')}
$$
\noindent is  an isomorphism of $F'$-vector-spaces. The right hand side is a tensor product of $\overline{\mathbf{Q}}$-vector spaces and the action of $\sigma\in\mathrm{Gal}(\overline{\mathbf{Q}}/F')$ is given by $\bigotimes_{s\in\{\tau_1,\dots,\tau_r\}}(x_{s}\otimes_{F,s} a_{s})\mapsto \bigotimes_{s\in\{\tau_1,\dots,\tau_r\}}\left(x_{s}\otimes_{F,\sigma s}\sigma(a_{s})\right)$.

\paragraph{Comparison isomorphisms} Let $I=\bigotimes_{j=1}^r I_{\tau_j}$, where 
$$I_{\tau_j}:M_{\mathrm{B},\tau_j}\otimes_{\mathbf{Q}}\mathbf{C}\overset{\sim}{\longrightarrow} M_{\mathrm{dR}}\otimes_{F,\tau_j}\mathbf{C}$$
\noindent are isomorphisms of $\mathbf{C}$-vector spaces, and $I'$ be the following isomorphism over $\mathbf{C}$ :
$$I':M_{\mathrm{B}}'\otimes_{\mathbf{Q}}\mathbf{C}\overset{\sim}{\longrightarrow} M_{\mathrm{dR}}'\otimes_{F'}\mathbf{C}.$$
\noindent The maps $I\circ\left(\mathscr I\otimes_{\mathbf{Q}}\text{id}_{\mathbf{C}}\right)$
and
$\left(\mathscr J\otimes_{F'}\text{id}_{\mathbf C}\right)\circ I'$ satisfy :

\begin{equation}
 \label{comparison}
\tag{$\star$}
I\circ\left(\mathscr I\otimes_{\mathbf Q}\text{id}_{\mathbf C}\right)=\left(\mathscr J\otimes_{F'}\text{id}_{\mathbf C}\right)\circ I' : M_{\mathrm B}'\otimes_{\mathbf Q}\mathbf C\overset{\sim}{\longrightarrow}\bigotimes_{j=1}^r\left(M_{\mathrm{dR}}\otimes_{F,\tau_j}\mathbf C\right).
\end{equation}

\noindent Yoshida's period conjecture consists of the isomorphisms $\mathscr I$, $\mathscr J$, $I$ and $I'$ satisfying $(\star)$. It is the Hodge-de Rham realization of the motivic conjecture above.

\paragraph{Complex conjugation : } Let $c_{\tau_j}$ be the complex conjugation on $M_{\mathrm{B},\tau_j}$. We will need the following hypothesis, which allows us to compare $c_{\tau_j}$ with $t_j^*$ on $M_{\mathrm{dR}}'\otimes_{F'}\mathbf C$.

\begin{hyp}
\label{goodaction}
The action of $t_j^*$ on $M_{\mathrm{dR}}'\otimes_{F'}\mathbf C$ corresponds via the isomorphism
$$\left(\mathscr I\otimes_{\mathbf Q}\mathrm{id}_{\mathbf C}\right)\circ \left(I'\right)^{-1} :M_{\mathrm{dR}}'\otimes_{F'}\mathbf C\longrightarrow M_{\mathrm{B}}'\otimes_{\mathbf Q}\mathbf C\longrightarrow\left(\bigotimes_{k=1}^r M_{\mathrm{B},\tau_k}\right)\otimes_{\mathbf Q}\mathbf C,$$
\noindent to the action of $c_{\tau_j}$ on $M_{\mathrm{B},\tau_j}$.
\end{hyp}

\subsection{Lattices and periods} 

Fix some $\omega_{\varphi}\neq 0$ in $F^rM_{\mathrm{dR}}'$. By definition of $M'$, there exists a finite set of places $S$ of $F$ such that for $v\notin S$, $T_v\omega_{\varphi}=a_v(E)\omega_{\varphi}$. 

Let $\Omega_{E/F}$ be the sheaf of differentials on $E/F$. Fix $\eta\neq 0\in H^0(E,\Omega_{E/F})=F^1M_{\mathrm dR}$. For $j\in\{1,\dots,n\}$, let 
$$\eta_j=\eta\otimes_{F,\tau_j} 1\in H^0\left(E\otimes_{F,\tau_j}{\overline{\mathbf{Q}}},\Omega_{(E\otimes_{F,\tau_j}{\overline{\mathbf{Q}}})/{\overline{\mathbf{Q}}}}\right)=\left(F^1M_{\mathrm{dR}}\right)\otimes_{F,\tau_j}{\overline{\mathbf{Q}}}.$$

\noindent Then
$$\bigotimes_{j=1}^r \eta_j\in\left(\bigotimes_{j=1}^r\left(F^1M_{\mathrm{dR}}\otimes_{F,\tau_j}{\overline{\mathbf{Q}}}\right)\right)^{\mathrm{Gal}({\overline{\mathbf{Q}}}/F')}=\mathscr J\left(F^r M_{\mathrm{dR}}'\right)$$
\noindent and there exists $\alpha\in F'^{\times}$ such that
$$\mathscr J(\alpha\omega_{\varphi})=\eta_1\otimes\dots\otimes\eta_r.$$

Let $j\in\{1,\dots,r\}$ and $E_j=E \otimes_{F,\tau_j}\mathbf C$. We shall denote by $H_1(E_j,\mathbf Z)^{\pm}$ the eigenspaces of the complex conjugation action on $H_1(E_j,\mathbf Z)$. Then
$$\left\{\int_{\Upsilon}\eta_j,\ \Upsilon\in H_1(E_{j},\mathbf Z)^{\pm}\right\}=\mathbf Z\Omega_j^{\pm},$$
\noindent where $\Omega_j^+\in\mathbf R\smallsetminus\{0\}$ and
$\Omega_j^-\in i\mathbf R\smallsetminus\{0\}$ are determined up to a sign. We fix the signs by imposing, e.g., $\mathrm{Re}\left(\Omega_j^+\right)>0$ and $\mathrm{Im}\left(\Omega_j^-\right)>0$.

Fix a character $\beta:\{1\}\times\prod_{j=2}^r\{\pm 1\}\rightarrow\{\pm 1\}$, and write $\beta=\prod_{j=2}^r\beta_j$. We set
$$\omega_{\varphi}^{\beta}=\left(\sum_{\sigma\in\{1\}\times\prod_{j=2}^r\{\pm 1\}} \beta(\sigma)t_{\sigma}^*\right)\omega_{\varphi}=\prod_{j=2}^r\left(1+\beta_j(-1)t_{j}^*\right)\omega_{\varphi}$$
\noindent and
$$\Omega^{\beta}=\prod_{j=2}^r\Omega_j^{\beta_j(-1)}.$$

The following identities 
$$\left(\bigotimes_{j=1}^r M_{\mathrm{B},\tau_j}\right)\otimes_{\mathbf Q}\mathbf C=\bigotimes_{j=1}^r\mathrm{Hom}_{\mathbf Z}(H_1(E_{j},\mathbf Z),\mathbf C)=\mathrm{Hom}_{\mathbf Z}\left(\bigotimes_{j=1}^r H_1(E_{j},\mathbf Z),\mathbf C\right)$$
\noindent and Yoshida's conjecture show that the image of $\alpha\omega_{\varphi}^{\beta}$ under the map
$$\left(\mathscr I\otimes_{\mathbf Q}\text{id}_{\mathbf C}\right)\circ I'^{-1}=I^{-1}\circ\left(\mathscr J\otimes_{F'}\text{id}_{\mathbf C}\right):M_{\mathrm{dR}}'\otimes_{F'} \mathbf C\longrightarrow\left(\bigotimes_{j=1}^r M_{B,\tau_j}\right)\otimes_{\mathbf Q} \mathbf C$$
\noindent is identified with the linear form
\begin{equation}
\label{formelineaireercompliquee}
\left\{\begin{array}{lcl}
          \bigotimes_{j=1}^r H_1(E_j,\mathbf Z)&\longrightarrow&\mathbf C\\
\Upsilon_1\otimes\dots\otimes\Upsilon_r&\longmapsto&\int_{\Upsilon_1\otimes\dots\otimes\Upsilon_r}\bigotimes_{j=1}^r\left(1+\beta_j(-1)t_{j}^*\right)\eta_j
         \end{array}\right.
\end{equation}

\noindent Hypothesis \ref{goodaction} allows us to be more explicit. Let $\Upsilon_1\otimes\dots\otimes\Upsilon_r\in \bigotimes_{j=1}^r H_1(E_j,\mathbf Z)$, then
\begin{align*}
\int_{\Upsilon_1\otimes\dots\otimes\Upsilon_r}\bigotimes_{j=1}^r\left(1+\beta_j(-1)t_{j}^*\right)\eta_j&=\left(\int_{\Upsilon_1}\eta_1\right)\prod_{j=2}^r\int_{\Upsilon_j}(1+\beta_j(-1)t_j^*)\eta_j\\&=
\left(\int_{\Upsilon_1}\eta_1\right)\prod_{j=2}^r\int_{\Upsilon_j+\beta_j(-1)c_j\Upsilon_j}\eta_j.
\end{align*}
\noindent and the linear form (\ref{formelineaireercompliquee}) takes values in $\Lambda_1\Omega^{\beta}=(\mathbf Z\Omega_1^++\mathbf Z\Omega_1^-)\Omega^{\beta}$.

Under the dual isomorphism $\mathscr I^*$ of $\mathscr I$, the lattices
$$\bigotimes_{j=1}^r{}_{\mathbf Z} H_1(E_{j},\mathbf Z)\subset\bigotimes_{j=1}^r{}_{\mathbf Q} M_{B,\tau_j}^*\quad\mathrm{and}\quad \mathrm{Im}\left(H_r(\mathrm{Sh}_H(G/Z,X)(\mathbf C),\mathbf Z)\longrightarrow (M_B')^*\right)$$
\noindent are commensurable. Thus there exists $\xi\in\mathbf{Z}\smallsetminus\{0\}$ such that 
$$\xi\mathrm{Im}\left(H_r(\mathrm{Sh}_H(G/Z,X)(\mathbf C), \mathbf Z)\longrightarrow (M_B')^*\right)\subset\mathscr I^*\left(\bigotimes_{j=1}^r{}_{\mathbf Z} H_1(E_{j}, \mathbf Z)\right)
.$$
\noindent This proves the following proposition :

\begin{prop}
\label{constructiondureseauprop}
Under the hypothesis made in this section ($E$ is modular, the multiplicity one in Yoshida's motivic conjecture and \ref{goodaction}), there exist  $\alpha\in F'^{\times}$ and $\xi\in\mathbf Z\smallsetminus \{0\}$ such that
$$\forall\gamma\in H_r(\text{Sh}_H(G,X)(\mathbf C),\mathbf Z),\ \forall \beta:\prod_{j=2}^r\{\pm 1\}\rightarrow\{\pm 1\},\qquad \xi\int_{\gamma}\alpha\omega_{\varphi}^{\beta}\in\Lambda_1\Omega^{\beta}.$$
\end{prop}

\subsection{General case}
\label{sectionmultiplicite}

When $m_H(\pi)=\mathrm{dim}\ \pi_f^H(\varphi)> 1$ Yoshida's conjecture is the following

\begin{conj}
The motive  $H^r(\mathrm{Sh}_H(G,X))^{(E)}$ is isomorphic to $\left(\bigotimes\limits_{\{\tau_1,\dots,\tau_r\}}\mathrm{Res}_{F/F'} M\right)^{m_H(\pi)}$.
\end{conj}

In general the motive $H^r(\mathrm{Sh}_H(G,X))^{(E)}$ has rank $\neq 2^r$. We shall provide Betti and de Rham realizations of a submotive $M'\subset H^r(\mathrm{Sh}_H(G,X))^{(E)}$ of rank $2^r$ and an isomorphism $M'\overset{\sim}{\longrightarrow} \bigotimes_{\{\tau_1,\dots,\tau_r\}}\mathrm{Res}_{F/F'} M$.

We need $0\neq\omega_{\varphi}\in F^r H_{\mathrm{dR}}^r(\mathrm{Sh}_H(G/Z,X)/F')^{(E)}$ satisfying 

\paragraph{de Rham cohomology} The $F'$-vector space
$$M'_{\mathrm{dR}}:=\left(\bigoplus_{\sigma\in\{\pm 1\}^r}\mathbf C t_{\sigma}^*(\omega_{\varphi}\otimes 1)\right)\cap H_{\mathrm{dR}}^r(\mathrm{Sh}_H(G/Z,X)/F')^{(E)}$$
\noindent has dimension $2^r$.

Thus
$$F^rM_{\mathrm{dR}}':=M_{\mathrm{dR}}'\cap F^r H^r_{\mathrm{dR}}(\mathrm{Sh}_H(G/Z,X)/F')^{(E)}=F'\omega_{\varphi}.$$

\paragraph{Betti cohomology} Let
$$I': H_{\mathrm{B}}^r(\mathrm{Sh}_H(G/Z,X)(\mathbf C),\mathbf Q)^{(E)}\otimes_{\mathbf Q}\mathbf C\overset{\sim}{\longrightarrow}  H_{\mathrm{dR}}^r(\mathrm{Sh}_H(G/Z,X)/F')^{(E)}\otimes_{F'}\mathbf C.$$

The $\mathbf Q$-vector space
$$M_{\mathrm{B}}':=I'^{-1}(M_{\mathrm{dR}}'\otimes_{F'}\mathbf C)\cap H_{\mathrm{B}}^r(\mathrm{Sh}_H(G/Z,X)(\mathbf C),\mathbf Q)^{(E)}$$
\noindent has dimension $2^r$.

\begin{definition}
\label{definitionrationnel}
 An element $\omega_{\varphi}\in F^r H_{\mathrm{dR}}^r(\mathrm{Sh}_H(G/Z,X)/F')^{(E)}$  is said rational if it satisfies the equations above.
\end{definition}

\paragraph{Comparison isomorphisms} There exist isomorphisms
$$\mathscr I : M_{\mathrm{B}}'\overset{\sim}{\longrightarrow} \bigotimes_{j=1}^r M_{\mathrm{B},\tau_j},$$
$$\mathscr J : M_{\mathrm{dR}}'\overset{\sim}{\longrightarrow} \left(\bigotimes_{j=1}^r(M_{\mathrm{dR}}\otimes_{F,\tau_j}\overline{\mathbf Q})\right)^{\mathrm{Gal}(\overline{\mathbf Q}/F')},$$
\noindent and
$$I_{\tau_j} : M_{\mathrm{B},\tau_j}\otimes_{\mathbf Q}\mathbf C\overset{\sim}{\longrightarrow}  M_{\mathrm{dR}}\otimes_{F,\tau_j}\mathbf C.$$
\noindent Set  $I=\bigotimes_{j=1}^r I_{\tau_j}$. We have

\begin{equation}
\label{egalitedesisodecomparaison}
\tag{$\star$} I\circ\left(\mathscr I\otimes_{\mathbf Q}\mathrm{id}_{\mathbf C}\right)=\left(\mathscr J\otimes_{F'}\mathrm{id}_{\mathbf C}\right)\circ I' : M_{\mathrm B}'\otimes_{\mathbf Q}\mathbf C\overset{\sim}{\longrightarrow} \bigotimes_{j=1}^r\left(M_{\mathrm{dR}}\otimes_{F,\tau_j}\mathbf C\right).
\end{equation}

As in Proposition \ref{constructiondureseauprop} we have
\begin{prop}
\label{propositiondedsreseaux}
Let $\omega_{\varphi}\in F^r H_{\mathrm{dR}}^r(\mathrm{Sh}_H(G/Z,X)/F')^{(E)}$ be rational. If $E$ is modular and if Yoshida's conjecture is true, then there exist $\alpha\in F'^{\times}$ and $\xi\in\mathbf Z\smallsetminus\{0\}$ such that 
$$\forall\gamma\in H_r(\mathrm{Sh}_H(G,X)(\mathbf C),\mathbf Z),\ \forall \beta:\prod_{j=2}^r\{\pm 1\}\rightarrow\{\pm 1\},\qquad \xi\int_{\gamma}\alpha\omega_{\varphi}^{\beta}\in\Lambda_1\Omega^{\beta}.$$
\end{prop}

\paragraph{Example}
 Let $H_1,H_2\subset \widehat{B}^{\times}$ be compact open subgroups such that there exists $g\in\widehat{B}^{\times}$ satisfying $g^{-1} H_1 g\subset H_2$. Let $\omega_{\varphi_2}\in F^rH^r_{\mathrm{dR}}(\mathrm{Sh}_{H_2}(G/Z,X)/F')^{(E)}$ be rational. Let us explain a way to obtain $\omega_{\varphi_1}\in F^rH^r_{\mathrm{dR}}(\mathrm{Sh}_{H_1}(G/Z,X)/F')^{(E)}$ rational.
 
Let $$\mathrm{pr}:\mathrm{Sh}_{g^{-1} H_1 g}(G/Z,X)\longrightarrow\mathrm{Sh}_{H_2}(G/Z,X)$$
\noindent be the map given by
 $$[x,b]_{g^{-1} H_1 g}\longmapsto[x,b]_{H_2}$$
\noindent and $$[\cdot g]:\mathrm{Sh}_{H_1}(G/Z,X)\longrightarrow\mathrm{Sh}_{g^{-1} H_1 g}(G/Z,X)$$
\noindent by $$[x,b]_{H_1}\longmapsto[x,bg]_{g^{-1} H_1 g}.$$
\noindent Let $\mathrm{pr}_g:\mathrm{Sh}_{H_1}(G/Z,X)\rightarrow\mathrm{Sh}_{H_2}(G/Z,X)$ be the composition of $\mathrm{pr}$ with $[\cdot g]$.

Choose $\theta_g\in\mathbf Q$. Set
$$\omega_{\varphi_1}:=\sum_{\substack{g\in\widehat{B}^{\times}\\ \mathrm{s.t.}\ g^{-1}H_1 g\subset H_2}}\theta_g \ \mathrm{pr}_g^*(\omega_{\varphi_2}),$$
$$(M_1')_{\mathrm{dR}}=\left(\sum_g\theta_g\ \mathrm{pr}_g^*\right)(M_2')_{\mathrm{dR}}$$ and $$(M_1')_{\mathrm{B}}=\left(\sum_g\theta_g\ \mathrm{pr}_g^*\right)(M_2')_{\mathrm B}.$$

\begin{prop}
 \label{projectiondeformesrationnelles}
If $\omega_{\varphi_1}\neq 0$, then the map \hspace{-0.2cm}$\displaystyle \sum_{\substack{g\in\widehat{B}^{\times}\\ \mathrm{s.t.}\ g^{-1}H_1 g\subset H_2}}\hspace{-0.4cm}\theta_g\mathrm{pr}_g^*$ is injective  on \mbox{$\bigoplus_{\sigma\in\{\pm 1\}^r}\mathbf C t_{\sigma}^*(\omega_{\varphi_2}\otimes 1)$} and $\omega_{\varphi_1}\in F^rH^r_{\mathrm{dR}}(\mathrm{Sh}_{H_1}(G/Z,X)/F')^{(E)}$ is rational.

\end{prop}

\begin{proof}
Assume that $\omega=\sum_{\sigma\in\{\pm 1\}^r}\lambda_{\sigma}t_{\sigma}^*\omega_{\varphi_2}\in\bigoplus_{\sigma\in\{\pm 1\}^r}\mathbf C t_{\sigma}^*(\omega_{\varphi_2}\otimes 1)$ (where $\lambda_{\sigma}\in\mathbf C$) is such that $\sum_g\theta_g\mathrm{pr}_g^*(\omega)=0$. We have the following equalities :
\begin{align*}
 \sum_g \theta_g\mathrm{pr}_g^*\omega&=\sum_g\theta_g\mathrm{pr}_g^*\sum_{\sigma} \lambda_{\sigma} t_{\sigma}^* \omega_{\varphi_2}\\
&=\sum_{\sigma}\lambda_{\sigma}t_{\sigma}^*\sum_g\theta_g\mathrm{pr}_g^*\omega_{\varphi_2}\\
\sum_g \theta_g\mathrm{pr}_g^*\omega&=\sum_{\sigma}\lambda_{\sigma}t_{\sigma}^*\omega_{\varphi_1}.
\end{align*}
Thus
$$\sum_{\sigma}\lambda_{\sigma}t_{\sigma}^*\omega_{\varphi_1}=0\in\bigoplus_{\sigma\in\{\pm 1\}^r} \mathbf C t_{\sigma}^*\omega_{\varphi_1},$$
\noindent and
$$\forall\sigma\in\{\pm 1\}^r\qquad \lambda_{\sigma}t_{\sigma}^*\omega_{\varphi_1}=0.$$
Hence $\forall \sigma\in\{\pm 1\}^r\ \lambda_{\sigma}\in 0$. The map $\sum_{g\in\widehat{B}^{\times}\ \mathrm{s.t.}\ g^{-1}H_1 g\subset H_2}\theta_g\mathrm{pr}_g^*$ commutes with $T_v,\ v\notin S$ and is an isomorphism $\bigoplus \mathbf C t_{\sigma}^*\omega_{\varphi_2}\rightarrow\bigoplus \mathbf C t_{\sigma}^*\omega_{\varphi_1}$. Hence $\omega_{\varphi_1}\in\left(\bigoplus_{\sigma\in\{\pm 1\}^r}\mathbf C t_{\sigma}^*(\omega_{\varphi_1}\otimes 1)\right)\cap F^r H_{\mathrm{dR}}^r(\mathrm{Sh}_{H_1}(G/Z,X)/F')^{(E)}$ is rational.

\end{proof}





\section{Toric orbits} 
Let $K/F$ be a quadratic extension satisfying the following properties :
\begin{enumerate}
 \item The places $\tau_2,\dots,\tau_r$ of $F$ are split in $K$.
 \item The places $\tau_1,\tau_{r+1},\dots,\tau_d$ are ramified in $K$.
\end{enumerate}

Thanks to the Skolem-Noether theorem, there exists an $F$-embedding $q:K\hookrightarrow B$,  unique up to conjugacy. We will denote by $q_j$ (resp. $\widehat{q}, q_{\mathbf{A}}$)
the induced embedding $K\hookrightarrow B_{\tau_j}$ (resp. $ \widehat{K}\hookrightarrow \widehat{B}$, $K_{\mathbf A}\hookrightarrow B_{\mathbf A}$). For each place $v$ of $F$, set $K_v=K\otimes_F F_v$.

\subsection{Cycles on $X$}
Let $T=\mathrm{Res}_{K/\mathbf Q}(\mathbf{G}_m)/\mathrm{Res}_{F/\mathbf Q}(\mathbf{G}_m).$ Thanks to Hilbert's  Theorem $90$ we have
$$T(A)=(K\otimes_{\mathbf Q} A)^{\times}/(F\otimes_{\mathbf Q} A)^{\times}$$
\noindent  for every $\mathbf Q$-algebra $A$.

Fix an embedding $q:T\hookrightarrow G/Z(G)$. The group $T(\mathbf R)$ is identified with $\prod_{j=1}^d K_{\tau_j}^{\times}/F_{\tau_j}^{\times}$ which allows us to define 
$q_j:K_{\tau_j}^{\times}/F_{\tau_j}^{\times}\longrightarrow G_{j,\mathbf R}.$

Let $\pi_0(T(\mathbf R))$ be the set of connected components of $T(\mathbf R)$ and denote by
$T(\mathbf R)^{\circ}$ the component of the identity. Fix a multi-orientation on $T(\mathbf R)^{\circ}=\prod_{j=1}^d(K_{\tau_j}^{\times}/F_{\tau_j}^{\times})^{\circ}$ (i.e. an orientation of each factor $(K_{\tau_j}^{\times}/F_{\tau_j}^{\times})^{\circ}$) and remark that
$$\pi_0(T(\mathbf R))=T(\mathbf R)/T(\mathbf R)^{\circ}\simeq\prod_{j=2}^r\{\pm 1\}.$$
\noindent We will focus on the orbits in $X$ under the action of $q(T(\mathbf R)^{\circ})$ by conjugation.

\begin{prop}

Let $\mathscr T^{\circ}$ be an orbit of $q(T(\mathbf R)^{\circ})$ in $X$. Then $\mathscr T^{\circ}$ decomposes into a product of orbits in $X_j$ under $q_j(T(\mathbf R)^{\circ})$ and is multi-oriented.

\end{prop}

\begin{proof}
The first part of this assertion follows from the natural decomposition $X=X_1\times\dots X_r$. The orbit $\mathscr T^{\circ}$ decomposes into orbits under $q_j((K_{\tau_j}^{\times}/F_{\tau_j}^{\times})^{\circ})$. For $j=1$, $q_j((K_{\tau_j}^{\times}/F_{\tau_j}^{\times})^{\circ})\simeq \mathbf S^1$ or a point and the orientation does not change. For $j\in\{2,\dots,r\}$, $q_j((K_{\tau_j}^{\times}/F_{\tau_j}^{\times})^{\circ})\simeq\mathbf R_+^{\times}$. The action of $\mathbf R_+^{\times}$ on itself by multiplication does not change the orientation. Hence the multi-orientation induced on $\mathscr T^{\circ}$ by $T(\mathbf R)^{\circ}$ is well-defined.

\end{proof}

In the following sections we shall fix some $q(T(\mathbf{R})^{\circ})$-orbit $\mathscr T^{\circ}$ whose projection on $X_1$ is a point. 

\begin{prop}
$\mathscr T^{\circ}$ is a connected multi-oriented submanifold of real dimension $r-1$.
\end{prop}

\begin{proof}
Recall that $\mathscr T^{\circ}$ is decomposed as $\mathscr T^{\circ}=\{z_1\}\times\mathscr T_2\times \dots\times\mathscr T_r$.
Fix $x\in X$ such that $\mathscr T^{\circ}=q(T(\mathbf{R})^{\circ}\cdot x$. Then for $j\in\{2,\dots,r\}$ we have $\mathscr T_j=q_j((K_{\tau_j}^{\times}/F_{\tau_j}^{\times})^{\circ})\cdot\mathrm{pr}_j(x)$. The group
$q_j((K_{\tau_j}^{\times}/F_{\tau_j}^{\times})^{\circ})$ is naturally identified with $\mathbf R_+^{\times}$ and $\mathscr T_j$ is a connected oriented manifold of real dimension one.

\end{proof}

As a corollary, we have the following decomposition :
$$\mathscr T^{\circ}=\{z_1\}\times\gamma_2\times\dots\times\gamma_r,$$
\noindent when $z_1$ is one of the two fixed points in the action of $q_1(T(\mathbf R)^{\circ})$ on $X_1$ and $\gamma_j$ is an oriented connected submanifold of real dimension one in $X_j$.

When we use the identification of $X$ with $(\mathbf C\smallsetminus \mathbf R)^r$, the action of $T(\mathbf R)$ on $X$ by conjugation is an action of $\mathrm{PGL}_2(\mathbf R)$ on $(\mathbf C\smallsetminus \mathbf R)^r$ by homography. Let $z\in K\smallsetminus F$. For $j\in\{2,\dots,r\}$ the matrix $q_j(z)$ is hyperbolic with exactly two fixed points in $\mathbf P^1(\mathbf R)$, $z_j$ and $z_j'$. The manifold $\gamma_j$ is then a circle arc in the Poincaré upper half-plane joining $z_j$ to $z_j'$ (or a line if $z_j'=\infty$). Figure 1 gives some examples of what could the $\gamma_j$s be in the case of circle arcs.

\begin{figure}[!h]
\label{orbites}
\centering
\includegraphics[height=4cm]{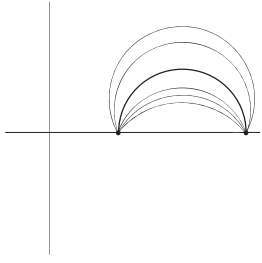} 
\caption{Case of circle arcs.}
\end{figure}

\subsection{Tori on $\mathrm{Sh}_H(G/Z,X)(\mathbf C)$}

Let $b\in\widehat{B}^{\times}$. We will denote by $\mathscr T_b^{\circ}$ the following subset of $\mathrm{Sh}_H(G/Z,X)(\mathbf C)$
$$\mathscr T_b^{\circ}=\left\{[x,b]_{H\widehat{F}^{\times}},x\in\mathscr T^{\circ}\right\}.$$

\begin{prop}
\label{Tbestuntore}
$\mathscr T_b^{\circ}$ is an oriented torus of real dimension $r-1$. 

\end{prop}

\begin{proof}
Let  $x,x'\in\mathscr T^{\circ}$ and $b\in\widehat{B}^{\times}$; we know that

\begin{align*}
 [x,b]_{H\widehat{F}^{\times}}=[x',b]_{H\widehat{F}^{\times}}&\Longleftrightarrow\exists k\in B^{\times}\ \mathrm{and}\ h\in H\widehat{F}^{\times}\qquad (kx',kbh)=(x,b)\\
&\Longleftrightarrow\exists k\in B^{\times}\cap b H\widehat{F}^{\times} b^{-1}\qquad kx'=x
\end{align*}
Since the projection of $\mathscr T^{\circ}$ on $X_1$ is a point, we have $k\in B\cap q_1(K_{\tau_1})=q_1(K)$ and
$$k\in q(K^{\times})\cap b H\widehat{F}^{\times} b^{-1}.$$ 
Thus the stabilizer $\mathscr W$ of $\mathscr T_b^{\circ}$ under the action of $q(K^{\times})$ is 
$$\mathscr W={q(K^{\times})\cap (b H\widehat{F}^{\times} b^{-1})}$$

\noindent which is commensurable with  $\mathcal O_{K,+}^{\times}/\mathcal O_F^{\times}$. This quotient has rank $r-1$ over $\mathbf Z$ as a consequence of  Dirichlet's units theorem :
$$\mathcal O_{K,+}^{\times}/\mathcal O_F^{\times}\simeq \mathrm{torsion}\times\mathbf Z^{r-1},$$
\noindent and the torsion is finite. The action of $T(\mathbf R)^{\circ}$ on $\mathscr T^{\circ}$ is given by $\prod_{j=2}^r(K_{\tau_j}^{\times}/F_{\tau_j}^{\times})^{\circ}$ and there is an isomorphism$$\prod_{j=2}^r(K_{\tau_j}^{\times}/F_{\tau_j}^{\times})^{\circ}\overset{\sim}{\longrightarrow}
\mathbf R^{r-1}.$$

\noindent The image $\widetilde{\mathcal O}$ of $\mathcal O_{K,+}^{\times}/\mathcal O_F^{\times}$ in $\mathbf R^{r-1}$ is isomorphic to $\mathbf Z^s$ with $s\leq r-1$. Denote by $\widetilde{\mathcal O}_K^{\times}$ the image of $\mathcal O_K^{\times}$ in $(K\otimes\mathbf R)^{\times,\ \mathrm{N}_{K/\mathbf Q}=1}$. As $$\prod_{j\notin\{2,\dots,r\}} K_{\tau_j}^{\times}/F_{\tau_j}^{\times}\qquad\mathrm{and}\qquad\frac{(K\otimes\mathbf R)^{\times,\ \mathrm{N}_{K/\mathbf Q}=1}}{\widetilde{\mathcal O}_K^{\times}}$$ \noindent are compact,
$\mathbf R^{r-1}/\widetilde{\mathcal O}$
\noindent is compact. Thus, the image of $\mathcal O_{K,+}^{\times}/\mathcal O_F^{\times}$ in $\mathbf R^{r-1}$ is a lattice.

The set $\mathscr T_b^{\circ}$ is a principal homogeneous space under $${q(K^{\times})}/{\mathscr W}\simeq (\mathbf R/\mathbf Z)^{r-1}.$$
\noindent It is a real torus in $\mathrm{Sh}_H(G/Z,X)(C)$ of dimension $r-1$, which is oriented by the fixed multi-orientation on $\mathscr T^{\circ}$.

\end{proof}

For each $u\in\pi_0(T(\mathbf R))$ and $b\in\widehat{B}^{\times}$ let 
$$\mathscr T_b^{u}=\left\{[q(u)\cdot x,b]_{H\widehat{F}^{\times}},\ x\in\mathscr T^{\circ}\right\}.$$
\noindent It is a real oriented torus of dimension $r-1$.

\begin{prop}
The set $$\{\mathscr T_b^{u}\mid\ b\in \widehat{B}^{\times},\ u\in\pi_0(T(\mathbf R))\}$$
\noindent does not depend on the choice of $q:K\hookrightarrow B$. 
\end{prop}

\begin{proof}
Let $\tilde{q}:K\hookrightarrow B$ be another embedding. Thanks to the Skolem-Noether theorem there exists $\alpha\in B^{\times}$ such that
 $$\forall k\in K\qquad \tilde{q}(k)=\alpha q(k)\alpha^{-1}.$$
\noindent Let $x_0\in X$, and assume that $\mathscr T^{\circ}=q(T(\mathbf R)^{\circ})\cdot x_0$. We have $\widetilde{\mathscr T}^{\circ}:=\tilde{q}(T(\mathbf R)^{\circ})\cdot \alpha(x_0)=\alpha\cdot\mathscr T^{\circ}$ and for each $u\in\pi_0(T(\mathbf R))$ 
$$\alpha\cdot q(u)\cdot\mathscr T^{\circ}=\tilde{q}(uT(\mathbf R)^{\circ})\cdot\alpha\cdot x_0.$$

Let $b\in\widehat{B}^{\times}$. As $\alpha\in B^{\times}$ we have 
$$\widetilde{\mathscr T}_b^{u}:=\left[\tilde{q}(u)\widetilde{\mathscr T}^{\circ},b\right]_{H\widehat{F}^{\times}}=\left[\alpha\cdot q(u)\cdot\mathscr T^{\circ},b\right]_{H\widehat{F}^{\times}}=\left[q(u)\cdot\mathscr T^{\circ},\alpha^{-1}\cdot b\right]_{H\widehat{F}^{\times}}=\mathscr T_{\alpha^{-1}b}^{u}.$$
\noindent The map $b\mapsto \alpha^{-1} b$ is a bijection. Thus 
$$\{\mathscr T_b^{u},\ b\in \widehat{B}^{\times},\ u\in\pi_0(T(\mathbf R))\}=\{\widetilde{\mathscr T}_b^{u},\ b\in \widehat{B}^{\times},\ u\in\pi_0(T(\mathbf R))\}.$$
\end{proof}

\subsubsection*{Action of \texorpdfstring{$\mathrm{Gal}(K^{\mathrm{ab}}/K)$}{Gal(Kab/K)}}

Let us denote by $K^{\mathrm{ab}}$ the maximal abelian extension of $K$ and by $\mathrm{rec}_K:K_{\mathbf A}^{\times}/K^{\times}\longrightarrow \mathrm{Gal}(K^{\mathrm{ab}}/K)$ the reciprocity map normalized by letting uniformizers correspond to geometric Frobenius elements.

The group $K_{\mathbf A}^{\times}$ acts on $\{\mathscr T_b^{u}\ |\ b\in \widehat{B}^{\times},\ u\in\pi_0(T(\mathbf R))\}$ by
$$\forall a=(a_{\infty},a_f)\in K_{\mathbf A}^{\times}=K_{\infty}^{\times}\times\widehat{K}^{\times}\ \forall b\in\widehat{B}^{\times}\qquad a\cdot\mathscr T_b^{u}=\mathscr T_{\widehat{q}(a_f)b}^{q(a_{\infty})u}.$$

\noindent The action of $k\in K^{\times}$ is trivial; as $q(k)\in B^{\times}$, the definition of $\mathrm{Sh}_H(G/Z,X)(\mathbf C)$ gives: 
$$k\cdot\mathscr T_b^u=[q(k)q(u)\mathscr T^{\circ},\widehat{q}(k)b]_{H\widehat{F}^{\times}}=[q(u)\mathscr T^{\circ},b]_{H\widehat{F}^{\times}}=\mathscr T_b^u.$$

\noindent The action of $F_{\mathbf A}^{\times}$ is trivial.
For $a=(a_{\infty},a_f)\in F_{\mathbf A}^{\times}$, and $b\in\widehat{B}^{\times}$, $\widehat{q}(a_f)b=b\widehat{q}(a_f)$ and $q(a_{\infty})q(u)\mathscr T^{\circ}=q(u)\mathscr T^{\circ}$ hence
$$a\cdot\mathscr T_b^u=[q(a_{\infty})q(u)\mathscr T^{\circ},\widehat{q}(a_f)b]_{H\widehat{F}^{\times}}=
[q(u)\mathscr T^{\circ},b]_{H\widehat{F}^{\times}}=\mathscr T_b^u.$$

\subsection{Special cycles on $\mathrm{Sh}_H(G/Z,X)(\mathbf C)$}

In this section we construct some $r$-chain on $\mathrm{Sh}_H(G/Z,X)(\mathbf C)$.

\begin{prop}
\label{cyclesdetorsionprop}
The homology class $[\mathscr T_b^{\circ}]\in H_{r-1}(\mathrm{Sh}_H(G/Z,X)(\mathbf C),\mathbf Z)$ of $\mathscr T_b^{\circ}$ is torsion. 
\end{prop}

\begin{proof}
Let us denote by $\mathrm{pr}$ the map
$$\mathrm{pr}:X\times\{b\}\longrightarrow\mathrm{Sh}_H(G/Z,X)(\mathbf{C}).$$
\noindent $\mathscr T_b^{\circ}$ is in the image of $\mathrm{pr}$ and 
$$\mathrm{pr}^{-1}(\mathscr T_b^{\circ})=(\{z_1\}\times\gamma_2\times\dots\times\gamma_r)\times\{b\}.$$

 Let $\omega\in H^{r-1}(\mathrm{Sh}_H(G/Z,X)(\mathbf{C}),\mathbf{C})$. Thanks to the Matsushima-Shimura theorem,  $\omega=\omega_{\mathrm{univ}}+\omega_{\mathrm{cusp}}.$ As $r-1\neq r$ we know that $\omega=\omega_{\mathrm{univ}}$.
\begin{itemize} 
\item  If $r-1$ is odd, then $H^{r-1}(\mathrm{Sh}_H(G/Z,X)(\mathbf{C}),\mathbf{C})=\{0\}$.
 
\item If $r-1=2s$ is even, $\omega$ is the pull-back of $\bigwedge_{j=2}^{r} \omega ^{(j)}$, where $$\omega^{(j)}= 1 \qquad\mathrm{or}\qquad
\frac{\mathrm{d} x_j\wedge\mathrm{d} y_j}{y_j^2}.$$ 
\noindent With the notations of the proof of Proposition \ref{Tbestuntore}, $\mathscr T_b^{\circ}$ is a principal homogeneous space under $\mathscr W$. Fix a fundamental domain $\widetilde{\mathscr W}$ of $\mathscr W$ in $\gamma_2\times\dots\times\gamma_r$. The incompatibility of degrees gives
 $$\int_{\mathscr T_b^{\circ}}\omega=\int_{\widetilde{\mathscr W}}\omega^{(2)}\wedge\dots\wedge\omega^{(r)}=0,$$
 \end{itemize} 
 
 $$\forall \omega\in H^{r-1}(\mathrm{Sh}_H(G/Z,X)(\mathbf{C}),\mathbf{C})\qquad \int_{\mathscr T_b^{\circ}}\omega=0.$$
\noindent This proves that $[\mathscr T_b^{\circ}]=0\in H_r(\mathrm{Sh}_H(G/Z,X)(\mathbf{C}),\mathbf{C})$ and $[\mathscr T_b^{\circ}]\in H_r(\mathrm{Sh}_H(G/Z,X)(\mathbf{C}),\mathbf{Z})$ is torsion.

\end{proof}

\begin{definition}
Let $n\in\mathbf Z_{>0}$ be the exponent of $H_{r-1}(\mathrm{Sh}_H(G/Z,X)(\mathbf C),\mathbf Z)_{\mathrm{tors}}$. Then
$$n[\mathscr T_b^{\circ}]=\partial \Delta_b^{\circ}$$
\noindent for some piece-wise differentiable $r$-chain $\Delta_b^{\circ}$.
\end{definition}

Proposition \ref{constructiondureseauprop} proves that the value of
$$\left(\frac{1}{\Omega^{\beta}}\xi\alpha\int_{\Delta_b^{\circ}}\omega_{\varphi}^{\beta}\right) \in\mathbf C$$
\noindent modulo $\Lambda_1$ does not depend on the particular choice of $\Delta_b^{\circ}$. If $T(\mathbf R)^{\circ}$ is fixed, then we have the following proposition.

\begin{prop}
\label{lemmedindependancedeschemins}
 Let $\mathscr T^{\circ}$ and $\mathscr T'^{\circ}$ be two special cycles such that $\mathrm{pr}_1(\mathscr T^{\circ})=\mathrm{pr}_1(\mathscr T'^{\circ})=\{z_1\}$. Assume that $\mathrm{pr}_j(\mathscr T^{\circ})$ and $\mathrm{pr}_j(\mathscr T'^{\circ})$ lie in the same connected component of $X_j$ for each $j\in \{2,\dots,r\}$. Let $n$ be the exponent of $H_{r-1}(\mathrm{Sh}_H(G/Z,X)(\mathbf C),\mathbf Z)_{\mathrm{tors}}$ and let $\Delta_b^{\circ}$ and $\Delta_b'^{\circ}$ satisfy
$$n[\mathscr T_b^{\circ}]=\partial\Delta_b^{\circ}\qquad\mathrm{and}\qquad n[\mathscr T_b'^{\circ}]=\partial \Delta_b'^{\circ}.$$
\noindent Then we have
$$\int_{\Delta_b^{\circ}}\omega_{\varphi}^{\beta}=\int_{\Delta_b'^{\circ}} \omega_{\varphi}^{\beta}\ 
(\mathrm{mod}\ \xi^{-1}\alpha^{-1}\Omega^{\beta}\Lambda_1).$$

\end{prop}

\begin{proof}
Our hypothesis allows us to decompose $\Delta_b'^{\circ}-\Delta_b^{\circ}$ into
$$\Delta_b'^{\circ}-\Delta_b^{\circ}=\mathrm{pr}(\{z_1\}\times \mathcal C) +\mathcal D,$$
\noindent where $\mathcal D$ is a cycle with $\partial \mathcal D=0$ and $\mathrm{pr}$ is the map 
$$\mathrm{pr}:\left\{\begin{array}{lcl}
X&\longrightarrow&\mathrm{Sh}_H(G/Z,X)(\mathbf C)\\
x&\longmapsto&[x,b]_{H\widehat{F}^{\times}}\end{array}
\right.$$

\begin{figure}[!h]
\hspace{-1.4cm}
\vspace{-1cm}
\centering
\includegraphics[width=16cm]{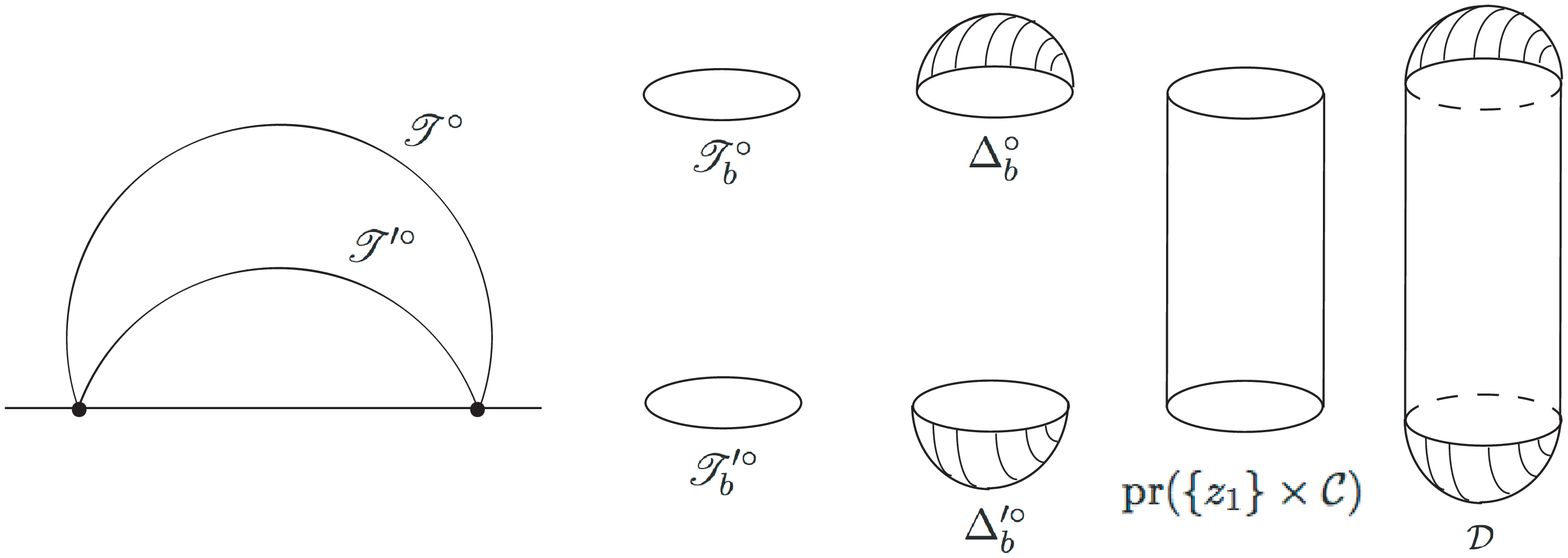}
\end{figure}

Let us show that $\int_{\Delta_b'^{\circ}-\Delta_b^{\circ}} \omega_{\varphi}^{\beta}\in \xi^{-1}\alpha^{-1}\Omega^{\beta}\Lambda_1$.
	
\noindent We have $$\omega_{\varphi}^{\beta}=\sum_{\varepsilon}\omega_{\varepsilon}\in\bigoplus_{\varepsilon:\{\tau_1,\dots,\tau_r\}\rightarrow\{\pm 1\}^r}\Gamma(\mathrm{Sh}_H(G/Z,X)(\mathbf C),(\Omega_H^{\mathrm{an}})^{\varepsilon}),$$
\noindent Each $\omega_{\varepsilon}\in\Gamma(\mathrm{Sh}_H(G/Z,X)(\mathbf C),(\Omega_H^{\mathrm{an}})^{\varepsilon})$ satisfies $$\mathrm{pr}^*(\omega_{\varepsilon})=\mathrm{d}z_1\wedge \omega_{\varepsilon}'$$ 
\noindent We have
$$\int_{\mathrm{pr}(\{z_1\}\times\mathcal C)}\omega_{\varepsilon}=\int_{\{z_1\}\times\mathcal C}\mathrm{d}z_1\wedge\omega_{\varepsilon}'=0,$$
\noindent thus
$$\int_{\{z_1\}\times\mathcal C}\omega_{\varphi}^{\beta}=0.$$

Thanks to Proposition \ref{constructiondureseauprop} we have
$$\int_{\mathcal D}\omega_{\varphi}^{\beta}\in \xi^{-1}\alpha^{-1}\Omega^{\beta}\Lambda_1$$
\noindent and the result follows.
\end{proof}

\begin{cor}
 The value modulo $\Lambda_1$ of
$$\left(\frac{1}{\Omega^{\beta}}\xi\alpha\int_{\Delta_b^{\circ}}\omega_{\varphi}^{\beta}\right) \in\mathbf C$$
\noindent  depends neither on the choice of $\mathscr T^{\circ}$ whose projection on $X_1$ is $\{z_1\}$ nor on $\Delta_b^{\circ}$ satisfying $n[\mathscr T_b^{\circ}]=\partial \Delta_b^{\circ}$.
\end{cor}

\begin{definition}
We set
$J_b^{\beta}=\frac{1}{\Omega^{\beta}}\xi\alpha\int_{\Delta_b^{\circ}}\omega_{\varphi}^{\beta}\ (\mathrm{mod}\ \Lambda_1) \in\mathbf C/\Lambda_1$, the image of $\mathscr T_b^{\circ}$ by an exotic Abel-Jacobi map.
\end{definition}

\subsubsection*{Properties of $J_b^{\beta}$}

For each $u\in\pi_0(T(\mathbf R))$ let $\Delta_b^u$ be some piece-wise differentiable chain satisfying $$n\left[[q(u)\cdot\mathscr T^{\circ},b]_{H\widehat{F}^{\times}}\right]=\partial \Delta_b^u.$$

\begin{prop}
We have
$$J_b^{\beta}=\frac{1}{\Omega^{\beta}}\xi\alpha\sum_{u\in\pi_0(T(\mathbf R))}\beta(u)\int_{\Delta_b^u}\omega_{\varphi}\ (\mathrm{mod}\ \Lambda_1).$$
\end{prop}

\begin{proof}
Let us identify $\pi_0(T(\mathbf R))$ with $\prod_{j=2}^r\{\pm 1\}$ and assume that the image of $T(\mathbf R)^{\circ}$ is $(1,\dots,1)$. Then
 $$\omega_{\varphi}^{\beta}=\sum_{u\in\pi_0(T(\mathbf R))} \beta(u) t_{u}^*(\omega_{\varphi}).$$
 \noindent The chains $t_u\Delta_b^{\circ}$ and $\Delta_b^u$ are in the same connected component. Thus using \ref{lemmedindependancedeschemins}, we have $$\int_{t_u\Delta_b^{\circ}}\omega_{\varphi}=\int_{\Delta_b^u}\omega_{\varphi}$$\noindent and the result follows.
 
 \end{proof}

Recall that $z_1\in X_1$ is fixed by $q(K_{\tau_1}^{\times})$. 

\begin{prop}
\label{changementdeTcirc}
 Let $\mathscr T^{\circ}$ and $\mathscr T'^{\circ}$ be two $q(T(\mathbf R)^{\circ})$-orbits such that $\mathrm{pr}_1(\mathscr T^{\circ})=\mathrm{pr}_1(\mathscr T'^{\circ})=\{z_1\}$. There exists a unique $u\in\pi_0(T(\mathbf R))$ such that, for all $j\in \{2,\dots,r\}$,
 $$\mathrm{pr}_j(\mathscr T'^{\circ})\ \mathrm{and}\ \mathrm{pr}_j(q(u)\cdot\mathscr T^{\circ})$$
 \noindent are in the same connected component of $X_j$.
 
If $J_b'^{\beta}\in\mathbf C/\Lambda_1$ denotes the value obtained from $\mathscr T'^{\circ}$, we have
$$J_b'^{\beta}=\beta(u)J_b^{\beta}.$$
\end{prop}

\begin{proof}
Let $x,x'\in X$ be such that $\mathscr T^{\circ}=q(T(\mathbf R)^{\circ})\cdot x$ (resp. $\mathscr T'^{\circ}=q(T(\mathbf R)^{\circ})\cdot x'$). There exists $u\in\pi_0(T(\mathbf R))$ such that for all $j\in \{1,\dots,r\}$, $\mathrm{pr}_j(q(u)\cdot x)$ and $\mathrm{pr}_j(x')$ are in the same connected component of $X_j$. As $\mathscr T'^{\circ}=q(u)\cdot\mathscr T^{\circ}$, the chain $\Delta_b'^{\circ}$ whose boundary up to torsion is $\left[\mathscr T'^{\circ},b\right]_{H\hat{F}^{\times}}$, equals $\Delta_b^u$. Thus
$$\sum_{u'\in\pi_0(T(\mathbf R))}\beta(u')\int_{\Delta_b'^{u'}}\omega_{\varphi}=\sum_{u'\in\pi_0(T(\mathbf R))}\beta(u')\int_{\Delta_b^{uu'}}\omega_{\varphi}=\beta(u)\sum_{u''\in\pi_0(T(\mathbf R))}\beta(u'')\int_{\Delta_b^{u''}}\omega_{\varphi}.$$
\end{proof}

 Let $q,q':K\hookrightarrow B$ be two embeddings and $x\in X$, $\mathscr T^{\circ}=q(T(\mathbf R)^{\circ})\cdot x$ (resp. $\mathscr T'^{\circ}=q'(T(\mathbf R)^{\circ})\cdot x'$). There exists $a\in B^{\times}$ such that $$q'=aqa^{-1}$$\noindent thanks to the Skolem-Noether theorem. For each $j\in\{1,\dots,r\}$, $\mathrm{pr}_j(\mathscr T^{\circ})$ and $\mathrm{pr}_j(\mathscr T'^{\circ})$ are in the same connected component of $X_j$ if and only if $\tau_j(\mathrm{nr}(a))>0$.
 
Using \ref{changementdeTcirc} we obtain

\begin{prop}
\label{corduchoixdeq}
If $$\alpha=(\mathrm{sgn}\circ\tau_j(\mathrm{nr}(a)))_{j\in\{1,\dots,r\}}\in\{\pm1\}^{r-1},$$\noindent then
$$J_b'^{\beta}=\beta(\alpha) J_b^{\beta}.$$

\end{prop}

Let $\mathrm{N}_{B^{\times}}(K^{\times})$ be the normalizer of $K^{\times}$ in $B^{\times}$. Let $a\in\mathrm{N}_{B^{\times}}(K^{\times})\smallsetminus K^{\times}$. After multiplying $a$ by an element in $K^{\times}$ we may assume
$$\forall j\in\{2,\dots,r\}\qquad \tau_j(\mathrm{nr}(a))>0.$$

\noindent We have
$$\mathrm{pr}_1(q(a)\cdot\mathscr T^{\circ})=t_1(z_1)$$
\noindent and
$$\forall j\in\{2,\dots,r\}\qquad\mathrm{pr}_j(q(a)\cdot\mathscr T^{\circ})=\mathrm{pr}_j(\mathscr T^{\circ})$$
\noindent but the orientations of $\mathrm{pr}_j(q(a)\cdot\mathscr T^{\circ})$ and $\mathrm{pr}_j(\mathscr T^{\circ})$ are not the same.

Thus
$$[t_1\mathscr T^{\circ},b]_{H\widehat{F}^{\times}}=[q(a)\mathscr T^{\circ},b]_{H\widehat{F}^{\times}}=[\mathscr T^{\circ},\widehat{q}(a)^{-1} b]_{H\widehat{F}^{\times}},$$
\noindent but the orientations differ by $(-1)^{r-1}$. Hence 

\begin{prop}
\label{actiondetun}
 The tori $\mathscr T_b^{\circ}$ and $t_1\mathscr T_{\widehat{q}(a) b}^{\circ}$ are the same up to orientation.
 \end{prop}





\section{Generalized Darmon's points}
\subsection{The main conjecture}

Let $\Phi_1:\mathbf C/\Lambda_1\overset{\sim}{\longrightarrow} E_1(\mathbf C)$
\noindent be the Weierstrass uniformization; i.e. the inverse of $\Phi_1$ is the Abel-Jacobi map for the differential $\eta_1$. For each $a_{\infty}\in K_{\infty}^{\times}$, fix some $r$-chain $q(a_{\infty})\cdot\Delta_b^{\beta}$ satisfying $n[q(a_{\infty})\cdot\mathscr T_b^{\beta}]=q(a_{\infty})\cdot\Delta_b^{\beta}$
\noindent and denote by $\beta(a_{\infty})$ the following sign 
$$\beta(a_{\infty})=\prod_{j=2}^r\beta\left(\mathrm{sgn}\left(\prod_{w\mid\tau_j}a_{\infty,w}\right)\right).$$

\begin{conj}
\label{conjectureprincipale}
The point
$$P_b^{\beta}=\Phi_1\left(\frac{1}{\Omega^{\beta}}\xi\alpha\int_{\Delta_b^{\beta}}\omega_{\varphi}\right)=\Phi_1(J_b^{\beta})\in E_1(\mathbf C)$$
\noindent lies in $E(K^{\mathrm{ab}})$ and
$$\forall a=(a_{\infty},a_f)\in K_{\mathbf A}^{\times}\qquad\mathrm{rec}_K(a)P_b^{\beta}=\Phi_1\left(\frac{\xi\alpha}{\Omega^{\beta}}\int_{q(a_{\infty})\cdot\Delta_{\widehat{q}(a_f)b}^{\beta}}\omega_{\varphi}\right)=\beta(a_{\infty}) P_{\widehat{q}(a_f)b}^{\beta}.$$
\end{conj}

\begin{remark}
The choice of $z_1\in X_1^{q_1(K_{\tau_1}^{\times})}$ fixes a morphism $h_1:\mathbf S\longrightarrow G_{1,\mathbf R},$
\noindent hence a morphism
$\mathbf C^{\times}=\mathbf S(\mathbf R)\longrightarrow G_{1,\mathbf R}(\mathbf R)=B_{\tau_1}^{\times}=(B\otimes_{F,\tau_1}\mathbf R)^{\times}$
\noindent satisfying $h_1(\mathbf C^{\times})=q_1(K_{\tau_1}^{\times})$. This fixes an embedding
$\tau_{1,K}:K\hookrightarrow\mathbf C$
\noindent such that the following diagram

$$
\xymatrix @!0 @R=20mm @C=3cm{
    \mathbf C^{\times} \ar@{->}[r]^{h_1} & (B\otimes_{F,\tau_1}\mathbf R)^{\times} \\
     & (K\otimes_{F,\tau_1}\mathbf R)^{\times} \ar@{->}[u]_{q_1} \ar@{->}[ul]^{\tau_{1,K}}
  }
$$

\noindent commutes. We may fix $\tilde{\tau}_1:K^{\mathrm{ab}}\hookrightarrow\mathbf C$
\noindent above $\tau_{1,K}$, such that 

$$
\xymatrix @!0 @R=20mm @C=3cm{
   F \ar@{->}[r]^{\tau_1}\ar@{->}[d] &\mathbf R \ar@{->}[r]&\mathbf C\\
   K\ar@{->}[urr]^{\tau_{1,K}} \ar@{->}[rr]& &K^{\mathrm{ab}}\ar@{.>}[u]_{\tilde{\tau}_1}
   }
$$

\noindent commutes. Moreover the isomorphism
$$\left\{\begin{array}{ccc}
\mathrm{Gal}(K^{\mathrm{ab}}/K)&\overset{\sim}{\longrightarrow}&\mathrm{Gal}(\tilde{\tau}_1(K^{\mathrm{ab}})/\tau_{1,K}(K))\\
\sigma&\longmapsto&\tilde{\tau}_1\circ\sigma\circ\tilde{\tau}_1^{-1}
\end{array}\right.$$
\noindent does not depend on the choice of $\tilde{\tau}_1$. If $\tilde{\tau}_1'$ is another embedding above $\tau_{1,K}$, then $\tilde{\tau}_1'=\tilde{\tau}_1\circ\sigma'$ with $\sigma'\in\mathrm{Gal}(K^{\mathrm{ab}}/K)$ and 
$$\forall\sigma\in\mathrm{Gal}(K^{\mathrm{ab}}/K)\qquad \tilde{\tau}_1'\circ\sigma\circ\tilde{\tau}_1^{\prime\ -1}=\tilde{\tau}_1\circ\sigma'\sigma\sigma^{\prime\ -1}\circ\tilde{\tau}_1^{-1}=\tilde{\tau}_1\circ\sigma\circ\tilde{\tau}_1^{-1}$$
\noindent because $\mathrm{Gal}(K^{\mathrm{ab}}/K)$ is commutative. Hence the Galois action of \ref{conjectureprincipale} does not depend on the particular choice of $\tilde{\tau}_1$.
\end{remark}

\begin{remark}
Using conjecture \ref{conjectureprincipale}, we obtain 

$$\forall a_{\infty}\in K_{\infty}^{\times}\qquad\mathrm{rec}_K(a_{\infty})P_b^{\beta}=\beta(a_{\infty})P_b^{\beta}.$$
$$\forall a\in F_{\mathbf A}^{\times}\qquad \mathrm{rec}_K(a)P_b^{\beta}=P_b^{\beta}.$$
\end{remark}

\subsection{Field of definition}

Let $B^{\times}_+=\{b\in B^{\times}\ |\ \forall j\in\{2,\dots,r\},\ \tau_j(\mathrm{nr}(b))>0\}$. It is diagonally embedded in $(B\otimes\mathbf R)^{\times}$. Set $$K_b^+=(K^{\mathrm{ab}})^{\mathrm{rec}_K(q_{\mathbf A}^{-1}(bH\hat{F}^{\times}b^{-1}B_+^{\times}))}\quad\mathrm{and}\quad K_b:=(K^{\mathrm{ab}})^{\mathrm{rec}_K(q_{\mathbf A}^{-1}(bH\hat{F}^{\times}b^{-1} B^{\times}))}\subset K_b^+.$$

\noindent Note that $K_b$ and $K_b^+$ depend on the choice of $q:K\hookrightarrow B$.

\begin{prop}
The point $P_b^{\beta}$ is defined over $K_b^+$ :  $P_b^{\beta}\in E(K_b^+)$.
\end{prop}

\begin{proof}
 Let $a=(1_{\infty},bhfb^{-1})(a_{\infty},1_f)\in q_{\mathbf A}^{-1}(bH\hat{F}^{\times}b^{-1} B_+^{\times})$ with $f\in\hat{F}^{\times}$ and $h\in H$. We have $$ \mathrm{rec}(a)P_b^{\beta}=\mathrm{rec}(q_{\mathbf A}^{-1}((1_{\infty},bhfb^{-1})) P_b^{\beta}
=P_{bhfb^{-1} b}^{\beta}=P_{bhf}^{\beta}
=P_b^{\beta}$$
\end{proof}

Remark that $\mathrm {rec}_K$ induces a surjection
 $$\mathcal R :\pi_0(T(\mathbf R))=\frac{(K\otimes_{\mathbf Q}\mathbf R)^{\times}}{(F\otimes_{\mathbf Q}\mathbf R)^{\times}(K\otimes_{\mathbf Q}\mathbf R)_+^{\times}}\simeq\prod_{j=2}^r\{\pm 1\}\twoheadrightarrow \mathrm{Gal}(K_b^+/K_b).$$
 
 \noindent Thus, we have
 
 \begin{prop}
 The points $P_b^{\beta}$ lie in $K_b^{\beta}=(K_b^+)^{\mathcal R(\mathrm{Ker}\ \beta)}$.
 \end{prop}

\begin{remark}
As $\mathrm{Ker} \beta$ has index 2 in $\prod_{j=2}^r\{\pm 1\}$, the field $K_b^{\beta}$ has degree 1 or 2 over $K_b$.
\end{remark}

Assume that the conductor $N$ of $E$ decomposes as $N=N_+N_-$ with $N=\mathfrak p_1\dots \mathfrak p_t$, $\mathfrak p_i$ distinct prime ideals of $\mathcal O_F$ and $t\equiv d-r\ \mathrm{mod}\ 2$. If $\mathrm{Ram}(B)=\{\tau_{r+1},\dots,\tau_d\}\cup\{\mathfrak p_1,\dots,\mathfrak p_t\}$ and $H=(R\otimes_{\mathbf Z}\widehat{\mathbf Z})^{\times}$ where $R\subset B$ is an Eichler order of level $N_+$, then $K_b$ is a ring class field of conductor $\mathfrak f_b$ and $K_b^+$ a ring class field of conductor $\mathfrak f_b\mathfrak f_{\infty}$, where $\mathfrak f_{\infty}=\prod_{j=2}^r \tau_j$.

\subsection{Local invariants of $B$}

Let $\pi$ be the irreductible automorphic representation of $B_{\mathbf A}^{\times}$  generated by $\varphi$ and
$$\eta_K=\eta_{K/F}:F_{\mathbf A}^{\times}/F^{\times}\mathrm{N}_{K/F}(K_{\mathbf A}^{\times})\longrightarrow\{\pm 1\}$$
\noindent the quadratic character of $K/F$. For each place $v$ of $F$ let $\mathrm{inv}_v(B_v)\in\{\pm 1\}$ be the invariant of $B$: $\mathrm{inv}_v(B_v)=1$ if and only if $B_v\simeq M_2(F_v)$.

Fix $b\in \widehat{B}^{\times}$ and  a character 
$$\chi:\mathrm{Gal}(K_b^+/K)\longrightarrow\mathbf C^{\times},$$
\noindent which will be identified with
$$K_{\mathbf A}^{\times}\overset{\mathrm{rec}_K}{\longrightarrow} \mathrm{Gal}(K^{\mathrm{ab}}/K)\longrightarrow\mathrm{Gal}(K_b^+/K)\overset{\chi}{\longrightarrow}\mathbf C^{\times}.$$

Let $L(\pi\times\chi,s)$ be the Rankin-Selberg $L$ function, see \cite{J} page 132 and \cite{JL} section 12. This function admits, since $\pi$ has trivial central character, 
 a holomorphic extension to $\mathbf C$ satisfying
 $$L(\pi\times\chi,s)=\varepsilon(\pi\times\chi,s)L(\pi\times\chi,1-s).$$

\noindent In this section, we prove the following

\begin{prop}
\label{determinationdeB}
 Let $b\in \widehat{B}^{\times}$ and assume conjecture \ref{conjectureprincipale}. If
$$e_{\overline{\chi}}(P_b^{\beta})=\sum_{\sigma\in\mathrm{Gal}(K_b^+/K)} \chi(\sigma)\otimes P_b^{\beta}\in E(K_b^+)\otimes \mathbf Z[\chi]$$
\noindent is not torsion, then $\beta=\chi_{\infty}$,
$$\forall v\neq\tau_1\qquad \eta_{K,v}(-1)\varepsilon(\pi_v\times\chi_v,\frac 12)=\mathrm{inv}_v(B_v)\quad\mathrm{and}\quad \varepsilon(\pi\times\chi,\frac 12)=-1.$$ 
\end{prop}

We shall use the following theorem (\cite{T} and  \cite{S}).

\begin{theorem}
\label{tunnelsaitowalds}
 The equality $\eta_{K,v}(-1)\varepsilon(\pi_v\times\chi_v,\frac 12)=\mathrm{inv}_v(B_v)$ holds if and only if there exists a non-zero invariant linear form
$$\ell_v:\pi_v\times\chi_v\longrightarrow\mathbf C$$
\noindent unique up to a scalar satisfying
 $$\forall a\in K_v^{\times}\ \forall u\in\pi_v\qquad \ell_v(q_v(a)u)=\chi_v(a)^{-1}\ell_v(u)$$
 \noindent i.e. $\ell_v$ is $q(K_v^{\times})$-invariant.

\end{theorem}

\begin{proof} (of Proposition \ref{determinationdeB}) 
We follow the proof of \cite{AN}, Proposition 2.6.2. 

Let $S'$ be a finite set of finite places of $F$ containing the places where $B$, $\pi$ or $K_b^+/F$ ramify, and such that the map $r=(r_v:K_v^{\times}\longrightarrow\mathrm{Gal}(K_b^+/K))_{v\in S'}$
\noindent obtained by composition
$$
r:\prod_{v\in S'} K_v^{\times}\longrightarrow K_{\mathbf A}^{\times}\overset{\mathrm{rec}_K}{\longrightarrow}\mathrm{Gal}(K^{\mathrm{ab}}/K)\longrightarrow \mathrm{Gal}(K_b^+/K)
$$
\noindent is surjective.

For each $v\in S'$ let
$$j_v:\left\{\begin{array}{lcl}
              K_v&\hookrightarrow&B_v\\
k&\longmapsto&b_v^{-1} q_v(k)b_v
             \end{array}\right.$$
\noindent and $$j=(j_v)_{v\in S'}:\prod_{v\in S'}K_v\hookrightarrow\prod_{v\in S'}B_v.$$

As $S'$ does not contain any archimedean place of $F$, 
$$\forall a\in\prod_{v\in S'} K_v^{\times}\qquad \left[\mathscr T^{\circ},\widehat{q}(a)b\right]_{H\hat{F}^{\times}}=\left[\mathscr T^{\circ},bj(a)\right]_{H\hat{F}^{\times}}$$
\noindent and
$$\forall a\in\prod_{v\in S'} K_v^{\times}\ \forall b\in \widehat{B}^{\times}\qquad \mathrm{rec}_K(a) P_b^{\beta}=P_{\widehat{q}(a)b}^{\beta}=P_{bj(a)}^{\beta}.$$

Let $(K_v^{\times})^{\circ}\subset K_v^{\times}$ be the inverse image of $(K_v^{\times}/\mathcal O_{K,v}^{\times})^{\mathrm{Gal}(K/F)}\subset K_v^{\times}/\mathcal O_{K,v}^{\times}.$

We have
$$K_v^{\times}/\mathcal O_{K,v}^{\times}F_v^{\times}\overset{\sim}{\longrightarrow}\left\{\begin{array}{ll}
                                                                  0&\mathrm{if}\ v\ \mathrm{is\ inert\ in}\ K/F\\
\mathbf Z/2\mathbf Z&\mathrm{if}\ v\ \mathrm{ramifies\ in}\ K/F\\
\mathbf Z&\mathrm{if}\ v\ \mathrm{splits\ in}\ K/F,
                                                                 \end{array}\right.$$
\noindent the quotient $(K_v^{\times})^{\circ}/F_v^{\times}$ is compact and
$$D_v:=K_v^{\times}/(K_v^{\times})^{\circ}\overset{\sim}{\longrightarrow}\left\{\begin{array}{ll}
                                                            \mathbf Z&\mathrm{if}\ v\ \mathrm{splits\ in}\ K/F\\
0&\mathrm{otherwise},
                                                           \end{array}\right.$$

$$(K_v^{\times})^{\circ}/\mathcal O_{K,v}^{\times}F_v^{\times}\overset{\sim}{\longrightarrow}\left\{\begin{array}{ll}
                                                                            \mathbf Z/2\mathbf Z&\mathrm{if}\ v\ \mathrm{ramifies\ in}\ K/F\\
0&\mathrm{otherwise.}
                                                                           \end{array}\right.$$

For each $v\in S'$, 
$C_v=\mathcal O_{K,v}^{\times}\cap\mathrm{Ker} (r_v)$
\noindent is  an open subgroup of $\mathcal O_{K,v}^{\times}$ and
$V_v^{\circ}= (K_v^{\times})^{\circ}/F_v^{\times}C_v$  is finite.

Let $V_v$ be the following subset of $K_v^{\times}/F_v^{\times}C_v$:
\begin{itemize}
 \item if $v$ does not split in $K/F$, 
$V_v^{\circ}=K_v^{\times}/F_v^{\times}C_v$
\noindent and $V_v:= V_v^{\circ}$. 

\item If $v$ splits in $K/F$, we fix some section of
$K_v^{\times}\twoheadrightarrow K_v^{\times}/(K_v^{\times})^{\circ}\overset{\sim}{\longrightarrow} \mathbf Z$.
Hence $K_v^{\times}=(K_v^{\times})^{\circ}\times D_v$
and there exists $n_v\geq 1$ such that
$\mathrm{Ker} (r_v\mid_{D_v})=n_v D_v.$

Fix a set of representatives $D_v'\subset D_v$ of $D_v/n_vD_v$ and set $V_v=V_v^{\circ}D_v'\subset K_v^{\times}/F_v^{\times}C_v.$
\end{itemize}

Let $V=\prod_{v\in S'} V_v\subset\prod_{v\in S'}K_v^{\times}/F_v^{\times}C_v$,
\noindent which is stable under multiplication by the abelian group $V^{\circ}=\prod_{v\in S'} V_v^{\circ}$ and such that
 $V\hookrightarrow \prod_{v\in S'}K_v^{\times}/F_v^{\times}C_v\overset{r}{\rightarrow} \mathrm{Gal}(K_b^+/K)$ is surjective with fibers of cardinality $\frac{|{V}|}{|{\mathrm{Gal}(K_b^+/K)|}}.$ We have
\begin{align*}
\frac{|{V}|}{|{\mathrm{Gal}(K_b^+/K)|}}e_{\overline{\chi}}(P_b^{\beta})&=\frac{|{V}|}{|{\mathrm{Gal}(K_b^+/K)|}}\sum_{\sigma\in \mathrm{Gal}(K_b^+/K)} \chi(\sigma)\otimes  \sigma\cdot P_b^{\beta}\\
&=\sum_{a\in V} \chi(a)\otimes P_{b j(a)}^{\beta}.
\end{align*}

Fix some open-compact subgroup $H_1\subset\bigcap_{a\in V} j(a) H j(a)^{-1}.$
Using the maps
$$\mathrm{Sh}_{H_1}(G/Z,X)\overset{[\cdot j(a)]}{\longrightarrow}\mathrm{Sh}_{j(a)^{-1}H_1j(a)}(G/Z,X)\overset{\mathrm{pr}}{\longrightarrow}\mathrm{Sh}_H(G/Z,X),$$
\noindent we have
\begin{align*}
\sum_{a\in V} \chi(a)\int_{\Delta_{bj(a)}^{\circ}}\omega_{\varphi}^{\beta}&=\sum_{a\in V} \chi(a)\int_{\Delta_{b}^{\circ}}[\cdot j(a)]^*\omega_{\varphi}^{\beta}\\
&=\int_{\Delta_{b}^{\circ}}\sum_{a\in V}\chi(a)[\cdot j(a)]^*\omega_{\varphi}^{\beta}\\
&=\int_{\Delta_{b}^{\circ}}\omega_1^{\beta},
\end{align*}
\noindent where $$\omega_1^{\beta}:=\sum_{a\in V}\chi(a)[\cdot j(a)]^*\omega_{\varphi}^{\beta}.$$

Whenever $\frac{|{V}|}{|{\mathrm{Gal}(K_b^+/K)|}} e_{\overline{\chi}}(P_b^{\beta})=\sum_{a\in V} \chi(a)\otimes P_{bj(a)}^{\beta}\in\mathbf Z[\chi]\otimes_{\mathbf Z}E(K_b^+)\subset \mathbf Z[\chi]\otimes_{\mathbf Z}\mathbf C/\Lambda_1$
is not torsion, there exists $\sigma:\mathbf Z[\chi]\hookrightarrow\mathbf C$ such that $$\frac{\xi\alpha}{\Omega^{\beta}}\int_{\Delta_b^{\circ}}\sum_{a\in V}^{} {}^{\sigma}\!\chi(a)[\cdot j(a)]^*\omega_{\varphi}^{\beta}\notin \mathbf Q[^{\sigma}\!\!\chi]\cdot\Lambda_1,$$\noindent where $^{\sigma}\!\chi=\sigma\circ\chi$.
The vector

$${}^{\sigma}\!\omega_1=\sum_{a\in V} {}^{\sigma}\!\chi(a)[\cdot j(a)]^*\omega_{\varphi}\in\pi^{H_1}\cap\Gamma(\mathrm{Sh}_{H_1}(G/Z,X),\Omega_{H_1})$$

\noindent is non-zero and invariant under $j(\prod_{v\in S'}(K_v^{\times})^{\circ})$. Moreover, $$\forall a\in\prod_{v\in S'}(K_v^{\times})^{\circ}\qquad j(a)\omega_1= {}^{\sigma}\!\chi^{-1}(a)\omega_1.$$

Let $$^{\sigma}\!\ell_{S'}:\bigotimes_{v\in S'}^{}{}^{\sigma}\!\pi_v=\bigotimes_{v\in S'} \pi_v\longrightarrow \mathbf C(^{\sigma}\!\chi^{-1})$$
\noindent be the $j(\prod_{v\in S'}(K_v^{\times})^{\circ})$-invariant projection on $\mathbf C\omega_1$.

Assume that $v\in S'$ does not split in $K$. In this case $(K_v^{\times})^{\circ}=K_v^{\times}$ and $^{\sigma}\!\ell_{S'}$ induces a  $q_v(K_v^{\times})$-invariant linear form ${}^{\sigma}\!\ell_v:\pi_v\rightarrow\mathbf C({}^{\sigma}\!\chi_v^{-1})$. We have $^{\sigma}\!\ell_v(\omega_{1,v})\neq 0$, where
$$\omega_{1,v}=\sum_{a_v\in V_v}^{} {}^{\sigma}\!\chi\circ r_v(a_v)[\cdot j_v(a_v)]^*\omega_{\varphi}.$$ 

As $\varepsilon_v(\pi_v\times ^{\sigma}\!\chi_v,\frac 12)$ is independent of $\sigma:\mathbf Z[\chi]\hookrightarrow\mathbf C$, Theorem \ref{tunnelsaitowalds} shows that
$$\eta_{K,v}(-1)\varepsilon(\pi_v\times\chi_v,\frac 12)=\mathrm{inv}_v(B_v).$$

\noindent When $v\in S'$ splits in $K$ or $v\notin S'\cup S_{\infty}$, the equality
$$\eta_{K,v}(-1)\varepsilon(\pi_v\times\chi_v,\frac 12)=1=\mathrm{inv}_v(B_v)$$
\noindent follows from calculations which may be found for example in \cite{N2} Proposition 12.6.2.4.

\paragraph{Global sign} If $v=\tau_j$ is an archimedean place, then $\varepsilon(\pi_v\times\chi_v,\frac 12)=1$. Moreover
$\eta_{K,v}(-1)=
 1\ \mathrm{if\ and\ only\ if}\ j\in\{2,\dots, r\}
$
\noindent and
$\mathrm{inv}_v(B_v)= 1\ \mathrm{if\ and\ only\ if}\ j\in\{1,\dots,r\}$.
Thus
$$\eta_{K,v}(-1)\mathrm{inv}_v(B_v)=\left\{\begin{array}{rl}
-1\times 1&\mathrm{if}\ j=1\\
1\times 1&\mathrm{if}\ j\in\{2,\dots,r\}\\
-1\times -1&\mathrm{}
\end{array}\right.$$
\noindent and
$$\forall j\in\{1,\dots,d\}\quad \varepsilon_v(\pi_v\times\chi_v,\frac 12)=\eta_{K,v}(-1)\mathrm{inv}_v(B_v)\times \left\{\begin{array}{rl}-1 & \mathrm{if}\ j=1 \\1 & \mathrm{if}\ j>1.\end{array}\right.$$
Hence
$$\varepsilon(\pi\times\chi,\frac 12)=-\prod_{v}\eta_{K,v}(-1)\mathrm{inv}_v(B_v)=-1.$$

\end{proof}

\subsection{Global invariant linear form and a conjectural Gross-Zagier formula}


For any open subgroup $H'\subset H$, $b\in\widehat{B}^{\times}$ and $u\in\pi_0(T(\mathbf R))$ fix $\Delta_{H',b}^u\in C^{r}(\mathrm{Sh}_H(G/Z,X)(\mathbf C),\mathbf Q)$ such that $\partial \Delta_{H',b}^u=[\mathscr T_{H',b}^u]$, where $\mathscr T_{H',b}^u=\{[q(u)x,b]_{H'\widehat{F}^{\times}},\ x\in\mathscr T^{\circ}\}$.

Recall that
$$\forall u'\in\pi_0(T(\mathbf R))\qquad t_{u'}\Delta_{H,b}^u=\Delta_b^{uu'}.$$

Let $\pi_{\infty}$ be the archimedean part of $\pi$. Fix $\varphi_{\infty}\in\pi_{\infty}$ a lowest weight vector of weight $(\underbrace{2,\dots,2}_{r},0,\dots,0)$ of $\pi_{\infty}$ and $\omega_{\varphi}$ such that
$\omega_{\varphi}=\varphi_{\infty}\otimes\varphi_f\in\pi_{\infty}\otimes\pi_f\subset S_2(B_{\mathbf A}^{\times}).$

Let us denote by ${}_{\mathbf Q}\pi_f$ the sub $\mathbf Q[\widehat{B}^{\times}]$-module of $\pi_f$ generated by $\varphi_f$.

\begin{prop}
The space ${}_{\mathbf Q}\pi_f$ is a $\mathbf Q$-vector space and
${}_{\mathbf Q}\pi_f\otimes_{\mathbf Q}\mathbf C\longrightarrow\pi_f$
\noindent is surjective.
\end{prop}

\begin{proof}
 The space $\mathrm{Im}({}_{\mathbf Q}\pi_f\otimes_{\mathbf Q}\mathbf C\rightarrow \pi_f)$ is a zero subvector space of $\pi_f$ invariant under $B_{\mathbf A}^{\times}$. As $\pi_f$ is irreducible, we have $\mathrm{Im}({}_{\mathbf Q}\pi_f\otimes_{\mathbf Q}\mathbf C\rightarrow \pi_f)=\pi_f$
 \noindent and ${}_{\mathbf Q}\pi_f\otimes_{\mathbf Q}\mathbf C\rightarrow \pi_f$ is surjective.
\end{proof}

Fix $\eta\neq 0\in H^0(E,\Omega_{E/F})$. There exists $\alpha\in F'^{\times}$ such that
$$\mathscr J(\alpha\omega_{\varphi})=\eta.$$
Fix a continuous character of finite order 
$\chi:K_{\mathbf A}^{\times}/K^{\times}F_{\mathbf A}^{\times}\longrightarrow\mathbf Z[\chi]^{\times}$. Let $H'\subset H$ be any open compact subgroup of $\widehat{B}^{\times}$
satisfying
$\chi\left(q_{\mathbf A}^{-1}(H' F_{\mathbf A}^{\times})\right)=1.$ Assume that there exists $b_0\in\widehat{B}^{\times}$ such that $b_0^{-1}H'b_0\subset H$. Let $\mathrm{pr}_{b_0}$ be the map $\mathrm{Sh}_{H'}(G/Z,X)\rightarrow\mathrm{Sh}_H(G/Z,X)$ defined on complex points~by
$$[x,b]_{H'\widehat{F}^{\times}}\mapsto [x,bb_0]_{H\widehat{F}^{\times}}.$$

\begin{prop}
\label{propabove}
If $b_0^{-1}H'b_0\subset H$ for some $b_0\in\widehat{B}^{\times}$,  then 
$$\forall Z'\in C^r(\mathrm{Sh}_{H'}(G/Z,X)(\mathbf C),\mathbf Z)\qquad \int_{Z'}\mathrm{pr}_{b_0}^*(\omega_{\varphi}^{\chi_{\infty}})\in\mathbf Q\alpha^{-1}\Omega^{\chi_{\infty}}\Lambda_1.$$

\end{prop}

\begin{proof}

Let $Z=\mathrm{pr}_{b_0}(Z')\in C^r(\mathrm{Sh}_{H}(G/Z,X)(\mathbf C),\mathbf Z)$. We have 
$$
\int_{Z'}\mathrm{pr}_{b_0}^*\omega_{\varphi}^{\chi_{\infty}}=\mathrm{deg}(\mathrm{pr}_{b_0}:Z'\rightarrow Z)\int_{Z}\omega_{\varphi}^{\chi_{\infty}}.
$$

\noindent Thanks to Proposition \ref{propositiondedsreseaux}, we have $\int_Z\omega_{\varphi}^{\chi_{\infty}}\in\mathbf Q\alpha^{-1}\Omega^{\chi_{\infty}}\Lambda_1$ hence $\int_{Z'}\mathrm{pr}_{b_0}^*\omega_{\varphi}^{\chi_{\infty}}\in\mathbf Q\alpha^{-1}\Omega^{\chi_{\infty}}\Lambda_1$.

\end{proof}

Denote by $\mathrm{pr}:\mathrm{Sh}_{H'}(G/Z,X)\longrightarrow\mathrm{Sh}_{H}(G/Z,X)$ the natural projection, and by $(K\otimes\mathbf R)^{\times}_+$ the set of elements in $(K\otimes\mathbf R)^{\times}$ whose  norm to $F$ is positive at each place of $F$. We have $\pi_0(T(\mathbf R))=\frac{(K\otimes\mathbf R)^{\times}}{(F\otimes\mathbf R)^{\times}(K\otimes\mathbf R)^{\times}_+}$.

\noindent  The following formula
$$
\ell_{\chi} (\omega') = \frac{1}{[H:H']\mathrm{deg}({\mathscr T}_{H',b} \overset{\mathrm{pr}}{\longrightarrow}
{\mathscr T}_{H,b})} \hspace{-0.2cm}\sum_{a\in \frac{K_{\mathbf A}^{\times}}{q_{\mathbf A}^{-1}(H'F_{\mathbf A}^{\times})(K\otimes \mathbf R)_+^{\times}}}\hspace{-0.8cm} \chi(a)\otimes
\int_{\Delta_{H',\widehat{q}(a_f)}^{q(a_{\infty})}} \omega'\quad (\mathrm{mod}\ \mathbf Q(\chi)\otimes_{\mathbf Q}\mathbf Q\alpha^{-1}\Omega^{\chi_{\infty}}\Lambda_1),
$$

\noindent where $\partial\Delta_{H',\widehat{q}(a_f)}^{q(a_{\infty})}=[\mathscr T_{H',\widehat{q}(a_f)}^{q(a_{\infty})}]$, is independent of the specific choice of $\Delta_{H',\widehat{q}(a_f)}^{q(a_{\infty})}$ : we can assume that $\omega'=\mathrm{pr}_{b_0}^*(\omega_{\varphi})$ for some $b_0\in\widehat{B}^{\times}$ ; decompose each $a\in K_{\mathbf A}^{\times}/q_{\mathbf A}^{-1}(H'F_{\mathbf A}^{\times})(K\otimes\mathbf R)_+^{\times}$ as $a=(a_f,1_{\infty})(1_f,a_{\infty})$. Remark that $$K_{\mathbf A}^{\times}/q_{\mathbf A}^{-1}(H'F_{\mathbf A}^{\times})(K\otimes\mathbf R)_+^{\times}=\widehat{K}^{\times}/\widehat{q}^{-1}(H'\widehat{F}^{\times})\times(K\otimes\mathbf R)^{\times}/(K\otimes\mathbf R)^{\times}_+,$$
\noindent  hence $a_f\in \widehat{K}^{\times}/\widehat{q}^{-1}(H'\widehat{F}^{\times})$ and $a_{\infty}\in (K\otimes\mathbf R)^{\times}/(K\otimes\mathbf R)^{\times}_+$.

\noindent Thanks to Proposition \ref{propabove}, the following formula 
\begin{align*}
\sum_{a_{\infty}\in K_{\infty}^{\times}}\chi_{\infty}(a_{\infty})\int_{\Delta_{H',\widehat{q}(a_f)}^{q(a_{\infty})}}\omega'&=\sum_{a_{\infty}\in K_{\infty}^{\times}}\chi_{\infty}(a_{\infty})\int_{\Delta_{H',\widehat{q}(a_f)}}\hspace{-0.4cm}t_{q(a_{\infty})}\mathrm{pr}_{b_0}^*\omega_{\varphi}\\
&=\int_{\Delta_{H,\widehat{q}(a_f)}}\hspace{-0.4cm}\omega_{\varphi}^{\chi_{\infty}}\qquad (\mathrm{mod}\ \mathbf Q\alpha^{-1}\Omega^{\chi_{\infty}}\Lambda_1)
\end{align*}
\noindent does not depend on the specific choice of $\Delta_{H',\widehat{q}(a_f)}^{q(a_{\infty})}$.

\noindent Thus, the expression of $\ell_{\chi}(\omega')$ above defines a linear form
$$
\ell_{\chi} : S_2^{H'} \cap \mathbf Q[\widehat{B}^{\times}] \omega_{\varphi} \longrightarrow
\mathbf Q(\chi) \otimes_{\mathbf Q} (\mathbf C/\mathbf Q\alpha^{-1}\Omega^{\chi_{\infty}}\Lambda_1).
$$

To simplify the notations, let
$$\delta_{H',H}=\mathrm{deg}({\mathscr T}_{H',b} \overset{\mathrm{pr}}{\longrightarrow}
{\mathscr T}_{H,b})
\qquad \mathrm{and}\qquad
W_{H'}=K_{\mathbf A}^{\times}/q_{\mathbf A}^{-1}(H'F_{\mathbf A}^{\times})(K\otimes \mathbf R)_+^{\times}.$$
Thus $$\ell_{\chi} (\omega') = \frac{1}{[H:H']\delta_{H',H}} \sum_{a\in W_{H'}} \chi(a)\otimes
\int_{\Delta_{H',\widehat{q}(a_f)}^{q(a_{\infty})}} \omega'.$$

\begin{prop}
 \begin{enumerate}
 \item Let $H''\subset H'\subset H$ be open compact subgroups such that $\chi(q_{\mathbf A}^{-1}(H'F_{\mathbf{A}}^{\times}))=~1$ and $\mathrm{pr}^*$ the map $\mathrm{pr}^*:S_2^{H'}(B_{\mathbf A}^{\times})\longrightarrow S_2^{H''}(B_{\mathbf A}^{\times}).$
 
 \noindent If $\omega'\in S_2^{H'}(B_{\mathbf A}^{\times})\cap\mathbf Q[\widehat{B}^{\times}]\omega_{\varphi}$, then
 $\ell_{\chi}(\omega')=\ell_{\chi}(\mathrm{pr}^*(\omega'))$
and $\ell_{\chi}$ defines a linear form on $\mathbf Q[\widehat{B}^{\times}]\omega_{\varphi}$.
\item We have
$$\forall a\in \widehat{K}^{\times}\ \forall \omega\in \mathbf Q[\widehat{B}^{\times}] \omega_{\varphi}\qquad \ell_{\chi}([\cdot \widehat{q}(a_f)]^* \omega)
= \chi_f(a)^{-1} \ell_{\chi}(\omega).$$
\item If $\chi$ factors through $\mathrm{Gal}(K_b^+/K)$ and if $P_b^{\beta}=\Phi_1\left(\int_{\Delta_{H,b}}\omega_{\varphi}^{\beta}\right)\otimes 1\in\mathbf C/\mathbf Q\Lambda_1$, then
$$e_{\overline\chi}
(P_b^{\chi_{\infty}})=\sum_{\mathrm{Gal}(K_b^+/K)}\chi(\sigma) \otimes \sigma (P_b^{\chi_{\infty}})\in\mathbf Q(\chi)\otimes_{\mathbf Q} E(K_b^+)\subset\mathbf Q(\chi)\otimes_{\mathbf Q}\left(\mathbf C/\mathbf Q\Lambda_1\right)$$
\noindent equals $\Phi_1(\ell_{\chi}([\cdot b]^* \omega_{\varphi}))$, up to a non-zero rational factor.
\end{enumerate}

\end{prop}

\begin{proof}
 \textbf{Proof of 1.} 
Let $a\in \widehat{K}^{\times}$. We have $\mathrm{pr}(\Delta_{H'',\widehat{q}(a_f)})=\Delta_{H',\widehat{q}(a_f)}$ and
$$\int_{\Delta_{H'',b}}\mathrm{pr}^*\omega'=\mathrm{deg}(\mathscr T_{H'',b}\longrightarrow \mathscr T_{H',b})\int_{\Delta_{H',b}}\omega'=\delta_{H'',H'}\int_{\Delta_{H',b}}\omega'.$$

As $\chi(q_{\mathbf A}^{-1}(H'F_{\mathbf A}^{\times}))=1$, we have (thanks to Proposition \ref{propabove})

\begin{align*}
 \ell_{\chi}(\mathrm{pr}^*\omega')&= \frac{1}{[H:H'']\delta_{H'',H}} \sum_{a\in W_{H''}} \chi(a)\otimes
\int_{\Delta_{H'',\widehat{q}(a_f)}^{q(a_{\infty})}} \mathrm{pr}^*\omega'& (\mathrm{mod}\ \mathbf Q(\chi)\otimes_{\mathbf Q}\mathbf Q\alpha^{-1}\Omega^{\chi_{\infty}}\Lambda_1)\\[4mm]
&= \frac{\delta_{H'',H'}}{\delta_{H'',H}} \sum_{a\in W_{H''}} \chi(a)\otimes
\int_{\Delta_{H',\widehat{q}(a_f)}^{q(a_{\infty})}} \omega'& (\mathrm{mod}\ \mathbf Q(\chi)\otimes_{\mathbf Q}\mathbf Q\alpha^{-1}\Omega^{\chi_{\infty}}\Lambda_1)\\[4mm]
&= \frac{\delta_{H'',H'}}{[H:H'']\delta_{H'',H}} \sum_{a\in W_{H'}} [H':H'']\chi(a) \otimes
\int_{\Delta_{H',\widehat{q}(a_f)}^{q(a_{\infty})}} \omega'& (\mathrm{mod}\ \mathbf Q(\chi)\otimes_{\mathbf Q}\mathbf Q\alpha^{-1}\Omega^{\chi_{\infty}}\Lambda_1)\\[4mm]
&= \frac{[H':H'']}{[H:H'']\delta_{H',H}} \sum_{a\in W_{H'}} \chi(a)\otimes
\int_{\Delta_{H',\widehat{q}(a_f)}^{q(a_{\infty})}} \omega'&(\mathrm{mod}\ \mathbf Q(\chi)\otimes_{\mathbf Q}\mathbf Q\alpha^{-1}\Omega^{\chi_{\infty}}\Lambda_1)\\[4mm]
&=\ell_{\chi}(\omega').
\end{align*}

\textbf{Proof of 2.}
Assume $H''$ is sufficiently small such that
$[\cdot \widehat{q}(a_f)]^*\mathrm{pr}^*\omega\in S_2^{H''}.$
 We have 

\begin{align*}
 \ell_{\chi} ([\cdot \widehat{q}(a_f)]^*\omega) &= \ell_{\chi} ([\cdot \widehat{q}(a_f)]^*\mathrm{pr}^*\omega)\\[4mm]
&= \frac{1}{[H:H'']\delta_{H'',H}} \sum_{a'\in W_{H''}} \chi(a')\otimes
\int_{\Delta_{H'',\widehat{q}(a')}^{q(a_{\infty}')}} [\cdot \widehat{q}(a_f)]^*\mathrm{pr}^* \omega& (\mathrm{mod}\ \mathbf Q(\chi)\otimes_{\mathbf Q}\mathbf Q\alpha^{-1}\Omega^{\chi_{\infty}}\Lambda_1)\\[4mm]
&= \frac{1}{[H:H'']\delta_{H'',H}} \sum_{a'\in W_{H''}} \chi(a')\otimes
\int_{\Delta_{H'',\widehat{q}(aa')}^{q(a_{\infty}')}}  \mathrm{pr}^*\omega& (\mathrm{mod}\ \mathbf Q(\chi)\otimes_{\mathbf Q}\mathbf Q\alpha^{-1}\Omega^{\chi_{\infty}}\Lambda_1)\\[4mm]
&= \frac{1}{[H:H'']\delta_{H'',H}} \sum_{a''\in W_{H''}} \chi(a''a^{-1})\otimes
\int_{\Delta_{H'',\widehat{q}(a'')}^{q(a_{\infty}'') }}\mathrm{pr}^*  \omega& (\mathrm{mod}\ \mathbf Q(\chi)\otimes_{\mathbf Q}\mathbf Q\alpha^{-1}\Omega^{\chi_{\infty}}\Lambda_1)\\[4mm]
&= \chi_f(a)^{-1}\frac{1}{[H:H'']\delta_{H'',H}} \sum_{a''\in W_{H''}} \chi(a'')\otimes
\int_{\Delta_{H'',\widehat{q}(a'')}^{q(a_{\infty}'') }} \mathrm{pr}^* \omega& (\mathrm{mod}\ \mathbf Q(\chi)\otimes_{\mathbf Q}\mathbf Q\alpha^{-1}\Omega^{\chi_{\infty}}\Lambda_1)\\[4mm]
&= \chi_f(a)^{-1} \ell_{\chi}(\mathrm{pr}^*\omega)\\[3mm]
&=\chi_f(a)^{-1} \ell_{\chi}(\omega)
\end{align*}

\textbf{Proof of 3.} 
As $\omega_{\varphi}\in S_2(B_{\mathbf A}^{\times})=\bigcup_H S_2^H(B_{\mathbf A}^{\times}),$
\noindent there exists $H'$ sufficiently small such that 
$$\omega_{\varphi}\in S_2^{H'}\quad\mathrm{and}\quad
[\cdot b]^*\omega_{\varphi}\in S_2^{H'}.$$
\noindent  Let $m=[K_{\mathbf A}^{\times}/q_{\mathbf A}^{-1}(H'F_{\mathbf A}^{\times})(K\otimes\mathbf R)_+^{\times}: \mathrm{Gal}(K_b^+/K)]$ and $\nu=\frac{1}{[H:H']\mathrm{deg}({\mathscr T}_{H'}\longrightarrow \mathscr T_{H})}.$ 
\noindent We have :

\begin{align*}
\ell_{\chi}(\circ[\cdot b]^* \omega_{\varphi})
&=\nu \sum_{a\in \frac{K_{\mathbf A}^{\times}}{q_{\mathbf A}^{-1}(HF_{\mathbf A}^{\times})(K\otimes \mathbf R)_+^{\times}}} \chi_f(a_f)\chi_{\infty}(a_{\infty})\otimes
\int_{\Delta_{H',\widehat{q}(a_f)}^{q(a_{\infty})}} [\cdot b]^*\omega_{\varphi}& (\mathrm{mod}\ \mathbf Q(\chi)\otimes_{\mathbf Q}\mathbf Q\alpha^{-1}\Omega^{\chi_{\infty}}\Lambda_1)\\[4mm]
&=\nu \sum_{a_f}\chi_f(a_f)\otimes\sum_{a_{\infty}} \chi_{\infty}(a_{\infty})
\mathrm{rec}_K(a_f)\cdot \int_{\Delta_{H',b}}\hspace{-0.3cm} t_{\mathrm{rec}_K(a_{\infty})}\omega_{\varphi}& (\mathrm{mod}\ \mathbf Q(\chi)\otimes_{\mathbf Q}\mathbf Q\alpha^{-1}\Omega^{\chi_{\infty}}\Lambda_1)\\[4mm]
&=\nu m\hspace{-0.4cm}\sum_{\sigma\in\mathrm{Gal}(K_b^+/K)}\chi(\sigma)\otimes\int_{\Delta_{H',b}}\sum_{a_{\infty}} \chi_{\infty}(a_{\infty})
 t_{\mathrm{rec}_K(a_{\infty})}\omega_{\varphi}& (\mathrm{mod}\ \mathbf Q(\chi)\otimes_{\mathbf Q}\mathbf Q\alpha^{-1}\Omega^{\chi_{\infty}}\Lambda_1)\\[4mm]
&=\nu m\hspace{-0.4cm}\sum_{\sigma\in\mathrm{Gal}(K_b^+/K)}\chi(\sigma)\otimes\int_{\Delta_{H',b}} \omega_{\varphi}^{\chi_{\infty}}& (\mathrm{mod}\ \mathbf Q(\chi)\otimes_{\mathbf Q}\mathbf Q\alpha^{-1}\Omega^{\chi_{\infty}}\Lambda_1),
\end{align*}

\noindent hence
$$e_{\overline{\chi}}(P_b^{\chi_{\infty}})=\Phi_1(\ell_{\chi}([\cdot b]^* \omega_{\varphi})).$$

\end{proof}

Let us consider the N\'eron-Tate height
$h_{\mathrm{NT}}:E(K^{\mathrm{ab}})\times E(K^{\mathrm{ab}})\longrightarrow\mathbf R$
extended to an hermitian form 
$$h_{\mathrm{NT}}:E(K^{\mathrm{ab}})\otimes\mathbf C\times E(K^{\mathrm{ab}})\otimes\mathbf C\longrightarrow\mathbf C.$$

\noindent Recall the condition 
\begin{equation}
\label{lestroispetitspoints}
\forall v\neq\tau_1\quad \varepsilon(\pi_v\times\chi_v,\frac 12)\eta_{K,v}(-1)=\mathrm{inv}_v(B)
\end{equation}

\noindent from Proposition \ref{tunnelsaitowalds}: if \ref{lestroispetitspoints} fails, then $P_b^{\chi_{\infty}}\in E(K^{\mathrm{ab}})$ is torsion.

\noindent In general, there should be some $k(b,\omega_{\varphi})\in\mathbf C$ such that
$$\forall\sigma : \mathbf Q(\chi)\hookrightarrow\mathbf C\quad h_{\mathrm{NT}}(e_{{}^{\sigma}\!\overline{\chi}}(P_b^{\chi_{\infty}}))=k(b,\omega_{\varphi}) L'(\pi\times {}^{\sigma}\!\chi,\frac 12),$$
\noindent as in Gross-Zagier, Zhang and Yuan-Zhang-Zhang \cite{GZ, Z1, YZZ}.

\noindent This formula explains the following conjecture :

\begin{conj}
Let  $K_{\chi}=(K^{\mathrm{ab}})^{\mathrm{Ker}(\chi)}$ be the extension of $K$ trivializing $\chi$.
 If $$\forall v\neq\tau_1\quad \varepsilon(\pi_v\times\chi_v,\frac 12)\eta_{K,v}(-1)=\mathrm{inv}_v(B),$$ \noindent then there exists $b\in\widehat{B}^{\times}$ such that $k(b,\omega_{\varphi})\neq 0$ and we have the following equivalences :
 
\begin{align*}
\ell_{\chi}\neq 0 &\Longleftrightarrow \exists b\in B_{\mathbf A}^{\times}\ \mathrm{such\ that}\ K_{\chi}\subset K_b^+\ \mathrm{and}\ e_{\overline{\chi}}(P_b^{\chi_{\infty}}) \in\mathbf Z[\chi]\otimes E(K_b^+)\ \mathrm{is\ not\ torsion}\\
&\Longleftrightarrow \exists \sigma:\mathbf Q(\chi)\hookrightarrow\mathbf C\qquad L'(\pi\times{}^{\sigma}\!\chi,\frac 12)\neq 0\\
&\Longleftrightarrow \forall \sigma:\mathbf Q(\chi)\hookrightarrow\mathbf C\qquad L'(\pi\times{}^{\sigma}\!\chi,\frac 12)\neq 0.
\end{align*}
\end{conj}

\section{A relation to Kudla's program}

The theorem of Gross-Kohnen-Zagier asserts that the positions of the traces to $\mathbf Q$ of classical Heegner points are given by the Fourier coefficients of some Jacobi form. The geometric proof of Zagier explained for example in \cite{Za} has been recently generalized by Yuan, Zhang and Zhang in \cite{YZZ} using a result of Kudla-Millson \cite{KM}. In this section we establish a relation between Darmon's construction and Kudla's program. This is a first step in an attempt to apply the arguments  of Zagier \cite{Za} and Yuan-Zhang and Zhang's \cite{YZZ} to Darmon's points.

\subsection{Some computations}

Let us fix a modular elliptic curve $E/F$ of conductor $N=N_+N_-$.
Assume $\mathrm{Ram}(B)=\{\tau_{r+1},\dots,\tau_d\}\cup\{v\mid N_-\}$
\noindent and that the quadratic extension $K/F$ satisfies the following hypothesis
$$\forall v\mid N_+\ \mathrm{splits\ in\ } K\qquad\forall v\mid N_-\ \mathrm{is\ inert\ in\ }K.$$

In particular, the relative discriminant $d_{K/F}$ is prime to $N$. Let $R$ be an Eichler order of $B$ of level $N_+$. Identify $K$ with its image in $B$ by $q$ and assume $K\cap R=\mathcal O_K$, 
\noindent $H=\hat{R}^{\times}$ (which implies that $\mathrm{dim}\pi_f^H=1$). 

Recall that $h_{z_1}$ defines an embedding $\tau_{1,K}:K\hookrightarrow\mathbf C$ and denote by $c$ the non-trivial element of $\mathrm{Gal}(K/F)$. Assume that Conjecture \ref{conjectureprincipale} is true for $\beta=1$ and let $P=\mathrm{Tr}_{K_1^+/K} P_1\in E(K).$

\begin{prop}
\label{calculdessignes}
 If $\varepsilon$ is the global sign of $E/F$, i.e. $\Lambda(E/F,s)=\varepsilon \Lambda(E/F,2-s)$, where $\Lambda$ is the completed $L$-function of $E/F$, then 
$c(P)=-\varepsilon P.$
\end{prop}

\begin{proof}
 Assume that $K=F(i)$ and $B=K(j)$, with $i^2=\mathfrak a\in F^{\times}$, $j^2=\mathfrak b\in F^{\times}$ and $ij=-ji$. Recall that
$$\mathscr T_1^{\circ}=[\mathscr T^{\circ},1]_{H\hat{F}^{\times}}$$
\noindent with $\mathscr T^{\circ}=\{z_1\}\times\gamma_2\times\dots\times\gamma_r$. Thus
$$c(\mathscr T_1^{\circ})=[\{t_1{z_1}\}\times\gamma_2\times\dots\times\gamma_r,1]_{H\hat{F}^{\times}}=(-1)^{r-1}[j^{-1}(\mathscr T^{\circ}),1]_{H\hat{F}^{\times}}$$
and $$c(\mathscr T_1^{\circ})=(-1)^{r-1}[\mathscr T^{\circ},j]_{H\hat{F}^{\times}}$$\noindent since $j\in B^{\times}$.
This shows that $c(P_1)=(-1)^{r-1}P_j$. We will write $P_j$ using only $P_1$. We will make the following abuse of language. For each place $v$ of $F$, $j_v$ shall denote the element $(1,\dots,1,\underbrace{j_v}_{v},1\dots)\in B_{\mathbf A}^{\times}$ and we will use the following lemma

\begin{lem}
\label{lelemme}
Let $b\in\widehat{B}^{\times}$ and $v$ a place of $F$. When $v\mid N_+$, set $k_v\in K_v^{\times}$ corresponding to $\begin{pmatrix}1&0\\0&\varpi_v^{\mathrm{ord}_v(N_+)}\end{pmatrix}$, where $\varpi_v$ is an uniformizer of $K_v$.  If $b_v=1$, then
$$P_{bj_v}=\left\{\begin{array}{ll}
-\varepsilon_vP_b&\mathrm{if}\ v\mid N_-\\
\varepsilon_v\mathrm{rec}_K(k_v^{-1})P_b&\mathrm{if}\ v\mid N_+\\
P_b&\mathrm{if}\ v\nmid N\end{array}\right.$$
\end{lem}

\begin{proof} (of the lemma)

For each $v$ inert in $K/F$ we have
\begin{align*}
 \mathrm{inv}_v(B)=1&\Longleftrightarrow B_v\simeq M_2(F_v)\\
&\Longleftrightarrow \mathfrak b\in\mathrm{N}_{K_v/F_v}(K_v^{\times})=\mathcal O_{F_v}^{\times}F_v^{\times}{}^2\\
&\Longleftrightarrow 2\mid\mathrm{ord}_v(\mathfrak b)
\end{align*}

\noindent As $\overline{j}=-j$, we have $\mathrm{nr}(j)=-j^2=-\mathfrak b$ and $$\mathrm{inv}_v(B)=1\Longleftrightarrow 2\mid\mathrm{ord}_v(\mathrm{nr}(j_v)).$$

\paragraph{If $v\mid N_-$} then $H_v=\mathcal O_{B_v}^{\times}$, where $\mathcal O_{B_v}$ is the unique  maximal order in $B_v$ hence $H_v\vartriangleleft B_v^{\times}$ and $B_v^{\times}/H_v^{\times}\simeq\mathbf Z$ by chosing some uniformizer. As $H_v$ is normal in $B_v^{\times}$, the map $$[\cdot j_v]:\mathrm{Sh}_H(G/Z,X)(\mathbf C)\longrightarrow\mathrm{Sh}_{j_v^{-1} Hj_v}(G/Z,X)(\mathbf C)$$
\noindent is well-defined on $\mathrm{Sh}_H(G/Z,X)(\mathbf C)$. Thus
$[\mathscr T^{\circ},bj_v]_{H\hat{F}^{\times}}=[\cdot j_v][\mathscr T^{\circ},b]_{H\hat{F}^{\times}}$ and
$$\int_{\Delta_{bj_v}^{\circ}}\omega_{\varphi}=\int_{\Delta_b^{\circ}}[\cdot j_v]^*\omega_{\varphi}=\int_{\Delta_b^{\circ}}\pi_v(j_v)\omega_{\varphi}.$$

Decompose $\pi=\pi(\varphi)=\otimes_v'\pi_v$. We have
$$\pi_v:B_v^{\times}\overset{\mathrm{nr}}{\longrightarrow}F_v^{\times}\overset{\mathrm{ord}_v}{\longrightarrow}\mathbf Z\longrightarrow\mathbf Z/2\mathbf Z\overset{\sim}{\longrightarrow}\{\pm 1\}.$$

\noindent Let us denote by $\alpha$ the following unramified character
$$\alpha : F_v^{\times}\overset{\mathrm{ord}_v}{\longrightarrow}\mathbf Z\longrightarrow\mathbf Z/2\mathbf Z\overset{\sim}{\longrightarrow}\{\pm 1\}$$\noindent satisfying $\pi_v=\alpha\circ\mathrm{nr}.$

\noindent As $v\mid N_-$,  $E$ has multiplicative reduction in $v$. The character $\alpha$ is trivial if and only if $E$ has split multiplicative reduction in $v$, i.e. $\varepsilon_v=-1$. 

Hence $$[\cdot j_v]^*\omega_{\varphi}=\alpha(\mathrm{nr}(j_v))\omega_{\varphi}=\left\{\begin{array}{ll}
 \omega_{\varphi}&\mathrm{if}\ \alpha=1\\
(-1)^{\mathrm{ord}_v(\mathrm{nr}(j))}\omega_{\varphi}&\mathrm{otherwise.}
\end{array}\right.$$

As $v\mid N_-$, $v\in\mathrm{Ram}(B)$ is inert in $K/F$ and $\mathrm{inv}_v(B)=-1$, thus $2\nmid\mathrm{ord}_v(\mathrm{nr}(j))$. Hence 
$$[\cdot j_v]^*\omega_{\varphi}=\alpha(\mathrm{nr}(j_v))\omega_{\varphi}=\left\{\begin{array}{ll}
 \omega_{\varphi}=-\varepsilon_v\omega_{\varphi}&\mathrm{if}\ \alpha=1\\
-\omega_{\varphi}=-\varepsilon_v\omega_{\varphi}&\mathrm{otherwise}
\end{array}\right.$$
\noindent and
$P_{bj_v}=-\varepsilon_v P_b.$

\paragraph{If $v\mid N_+$} 
then we fix some uniformizer $\varpi_v$ of $F_v$ and an isomorphism $B_v\simeq \mathrm{M}_2(F_v)$ which identifies $K_v$ with the set of diagonal matrices and $R_v$ with $\left\{\begin{pmatrix} a&b\\ c&d\end{pmatrix}\in M_2(\mathcal O_{F,v})\ \left|\ \varpi_v^{\mathrm{ord}_v(N_+)}\mid c\right.\right\}$.

As $\mathrm{inv}_v(B_v)=1$,  $j_v$ is a local norm. There exists $k_v\in K_v$ such that $j_v=\mathrm{N}_{K_v/F_v}(k_v)$. We may assume that $j_v^2=1$. Moreover $j_v$ is in the normalizer of $K_v^{\times}$ in $B_v^{\times}$ we thus identify $j_v$ to~$
\begin{pmatrix}
0&1\\1&0
   \end{pmatrix}$.

Set $$W_v=\begin{pmatrix}
              0&1\\ \varpi_v^{\mathrm{ord}_v(N_+)}&0
             \end{pmatrix}=\begin{pmatrix}
                1&0\\
0& \varpi_v^{\mathrm{ord}_v(N_+)}
                \end{pmatrix}
                \begin{pmatrix}
0&1\\1&0
   \end{pmatrix}=k_vj_v.
$$
\noindent This matrix is in the normalizer of $R_v$ in $B_v$. As $W_v$ normalize $H_v$, 
$$[\mathscr T^{\circ},bj_v]_{H\hat{F}^{\times}}=[\mathscr T^{\circ},bk_v^{-1}W_v]_{H\hat{F}^{\times}}=[\cdot W_v][\mathscr T^{\circ},bk_v^{-1}]_{H\hat{F}^{\times}}.$$

Decompose $\omega_{\varphi}=\bigotimes_{v\mid N_+}\omega_v\otimes\omega'$, where  $\omega_v$ satisfies$[\cdot W_v]^*\omega_v=\varepsilon_v\omega_v$; then 
$$\int_{\Delta_{bj_v}^{\circ}}\omega_{\varphi}=\varepsilon_v\int_{\Delta_{bk_v^{-1}}^{\circ}}\omega_{\varphi}.$$
\noindent As $b_v=1$, $$P_{bj_v}=\varepsilon_v\mathrm{rec}_K(k_v^{-1})P_b.$$

 \paragraph{If $v\nmid N$} then by a similar calculation 
  we obtain $$P_{bj_v}=\mathrm{rec}_K(k_v^{-1})P_b.$$

\end{proof}

\paragraph{End of the proof of Proposition \ref{calculdessignes}} Lemma \ref{lelemme} implies that  $$c(P_1)=(-1)^{r-1}\prod_{v\mid N_-} (-\varepsilon_v)\prod_{v\mid N_+} \varepsilon_v\mathrm{rec}_K(k_v^{-1})P_1$$
\noindent and 
$$\forall a\in K_{\mathbf A}^{\times}\qquad c(\mathrm{rec}_K(a)P_1)=(-1)^{r-1}\prod_{v\mid N_-} (-\varepsilon_v)\prod_{v\mid N_+} \varepsilon_v\mathrm{rec}_K(k_v^{-1})\mathrm{rec}_K(a)P_1.$$

\noindent As $P\in E(K)$, we know that $\mathrm{rec}_K(k^{-1})P=P$. Thus
\begin{equation}
\label{signedelequation}
c(P)=(-1)^{r-1}\prod_{v\mid N_-} (-\varepsilon_v)\prod_{v\mid N_+} \varepsilon_v P=(-1)^{r-1}(-1)^{|{\{v\mid N_-\}|}}\prod_{v\nmid\infty}\varepsilon_v P.
\end{equation}

We have to show that $(-1)^{r-1}\prod_{v\mid N_-} (-\varepsilon_v)\prod_{v\mid N_+} \varepsilon_v =-\varepsilon$. For each $v\mid\infty$ we have $\varepsilon_v=-1$. Since $\prod_{v\mid\infty}=(-1)^d$, the sign in equation (\ref{signedelequation}) is $$(-1)^d\underbrace{\prod_v \varepsilon_v}_{=\varepsilon} (-1)^{r-1}(-1)^{|{\{v\mid N_-\}|}}.$$

Recall that $\{v\mid N_-\}=\mathrm{Ram}(B)\cap S_f$. As $|\mathrm{Ram}(B)|$ is even, we have
$$(-1)^{|{\{v\mid N_-\}|}}=(-1)^{|{\mathrm{Ram}(B)\cap S_{\infty}|}}=(-1)^{d-r}.$$
\noindent Hence

$$c(P)=(-1)^d\varepsilon (-1)^{r-1}(-1)^{|{\{v\mid N_-\}|}} P=-\varepsilon P.$$

\end{proof}

\begin{remark}

The above computations are a particular case of a result of Prasad, \cite{Prasad} Theorem 4, which asserts that if $\mathrm{Hom}_{K_v^{\times}}(\pi_v,\mathbf{1})\neq\{0\}$, then the non trivial element in $\mathrm{N}_{B_v^{\times}}(K_v^{\times})\backslash K_v^{\times}$ acts on $\mathrm{Hom}_{K_v^{\times}}(\pi_v,\mathbf{1})$ by multiplication by $\mathrm{inv}_v(B)\varepsilon_v=\mathrm{inv}_v(B)\varepsilon(\pi_v,\frac 12)\in\{\pm 1\}$.

\end{remark}

\subsection{Orthogonal Shimura manifolds}

Until the end of this paper we shall assume $h_F^+=1$.

Let us recall some definitions used by Kudla \cite{K} in the particular case $r=1$.
Let $n\in\mathbf Z_{\geq 0}$ and let $(V,Q)$ be a quadratic space over $F$ of dimension $n+2$. We assume that the signature of $V\otimes_{F,\tau_j}\mathbf R$ is
$$(n,2)\times(n+1,1)^{r-1}\times(n+2,0)^{d-r}.$$
Denote by $D$ the symmetric space of $G=\mathrm{Res}_{F/\mathbf Q}\mathrm{GSpin}(V)$. $D$ is the product of the oriented symmetric spaces of $V_j=V\otimes_{\tau_j,F}\mathbf R$. Thus 
$D=D_1\times\dots D_d$,
\noindent where $D_j$ is the set of oriented positive subspaces in $V_j$ of maximal dimension.
 For each $x\in V$ let $x_j$ be the image of $x$ in $V_j$. Assume that $Q(x)$ is totally positive.
 Set $V_x=x^{\perp}$, $G_x=\mathrm{Res}_{F/\mathbf Q}\mathrm{GSpin}(V_x)$ and for each $j\in\{1,\dots,d\}$
$$D_{x_j}=\{z\in D_j\ z\perp x_j\}.$$

We shall focus on the following real cycle on the Shimura manifold $G(\mathbf Q)\backslash D\times G(\widehat{\mathbf Q})/H$.

\begin{definition}
 Let $H$ be an open compact subgroup in $G(\widehat{\mathbf Q})$ and $g\in G(\widehat{\mathbf Q})$. The cycle $Z(x,g;H)$ is defined to be the image of the map
$$Z(x,g;H):\left\{\begin{array}{ccc}
 G_x(\mathbf Q)\backslash D_x\times G_x(\widehat{\mathbf Q})/H_x^g&\longrightarrow&G(\mathbf Q)\backslash D\times G(\widehat{\mathbf Q})/H\\
G_x(\mathbf Q)(y,u)H_x^g&\longmapsto&G(\mathbf Q)(y,ug)H\widehat{F}^{\times},
\end{array}\right. $$
\noindent where $H_x^g$ denotes $G_x(\widehat{\mathbf Q})\cap g Hg^{-1}$.
\end{definition}

\paragraph{Example (including Proposition \ref{propositionexemple}) :}
Fix $D_0\in F$ satisfying
$$\left\{\begin{array}{ll}
  \tau_j(D_0)>0&\mathrm{if}\ j\in\{1,r+1,\dots,d\}\\
 \tau_j(D_0)<0&\mathrm{if}\ j\in\{2,\dots,r\}
  \end{array}\right.$$

\noindent Set $$(V,Q)=(B^{\mathrm{Tr=0}}, D_0\cdot \mathrm{nr}).$$
\noindent 

\noindent $(V\otimes_{F,\tau_j}\mathbf R,\tau_j\circ D_0\cdot\mathrm{nr})$ has signature
$$\left\{\begin{array}{ll}
(1,2)&\mathrm{if}\ j=1\\
(2,1)&\mathrm{if}\ j\in\{2,\dots,r\}\\
(3,0)&\mathrm{if}\ j\in\{{r+1},\dots,{d}\}.\end{array}\right.$$

 Let $G=\mathrm{Res}_{F/\mathbf Q}\mathrm{GSpin}(V)$. The action of $B^{\times}$ on $V$ by conjugation induces an isomorphism
$$\begin{array}{lcl}B^{\times}&\overset{\sim}{\longrightarrow}&\mathrm{GSpin}(V)\\
   b&\longmapsto&(v\mapsto bvb^{-1}),
  \end{array}
$$
\noindent thus $G\simeq\mathrm{Res}_{F/\mathbf Q}(B^{\times})$.

Let $x\in V$ such that
$Q(x)\gg 0$, and denote by $x_j$ its image in $V\otimes_{F,\tau_j}\mathbf R$. Denote by $K$ the quadratic extension $F+Fx$ and $T=\mathrm{Res}_{K/\mathbf Q}(\mathbf{G}_m)/\mathrm{Res}_{F/\mathbf Q}(\mathbf{G}_m)$ as above. Let $q$ be the inclusion $K\hookrightarrow V\rightarrow B$.

\begin{prop}
\label{propositionexemple}
The set
$$D_x=D_{x_1}\times\dots\times D_{x_r}$$
\noindent is a $q(T(\mathbf R))^{\circ}$-orbit in $D$ whose projection on $D_1$ is a point. \end{prop}

\begin{proof}
As $x\in V$,  $\mathrm{Tr}(x)=0$ and $x^2=-\mathrm{nr}(x)=-\frac{Q(x)}{D_0}\in F^{\times}$. Let $j\in\{1,\dots,r\}$. We have $\tau_j(Q(x))>0$ hence $\tau_j(x^2)\tau_j(D_0)<0$. Thus $\tau_1$ ramifies in $K$ and $\tau_2,\dots,\tau_r$ are split. Moreover $q_1(K^{\times})$ fixes $x_1$ by defintion of $K$.

\end{proof}


Let us focus on the general case when $V$ has dimension $n$.
Fix $t\in F$ satisfying $\forall j\in\{1,\dots,r\}\ \tau_j(t)>0.$
\noindent $G(\widehat{\mathbf Q})$ acts on $\Omega_t=\{x\in V(F)\ |\quad Q(x)=t\}$ by conjugation.

 Let $\varphi$ be a Schwartz function on $V(\widehat{F})$. Assume $\Omega_t\neq\emptyset$ and fix $x\in\Omega_t$. Denote by $Z(y,\varphi; H)$ the following sum 
$$Z(t,\varphi;H)=\sum_{g\in G_x(\widehat{\mathbf Q})\backslash G(\widehat{\mathbf Q})/H\widehat{F}^{\times}}\varphi(g^{-1}\cdot x)Z(x,g;H).$$

\noindent Proposition \ref{cyclesdetorsionprop} showed that for $n=1$ $[Z(x,g;H)]=0\in H_{r-1}(\mathbf{Sh}_H(G/Z,X)(\mathbf C),\mathbf C)$. A natural invariant to consider is the refined class
$$\{Z(t,\varphi;H)\}=\omega\mapsto J_b^{\beta}\in\frac{(\mathrm{Harm}^r(\mathrm{Sh}_{H}(G/Z,X)(\mathbf C))^*}{\mathrm{Im}(H_r(\mathrm{Sh}_H(G/Z,X)(\mathbf C),\mathbf Z)\rightarrow\mathrm{Harm}^r(\mathrm{Sh}_H(G/Z,X)(\mathbf C))^*)},$$ 
\noindent where $\mathrm{Harm}^r(\mathrm{Sh}_H(G/Z,X)(\mathbf C))$ is the set of harmonic differential forms on $\mathrm{Sh}_H(G/Z,X)(\mathbf C)$.

 In order to adapt the work of Yuan, Zhang and Zhang, we need the following conjecture

\begin{conj}
\label{conjecturesixdeuxtrois}
In the situation of the above example $(V,Q)=(B^{\mathrm{Tr}=0},D_0\cdot\mathrm{nr})$, the sum
$$\sum_{\substack{t\in\mathcal O_F\\ t\gg 0}}\{Z(t,\varphi;H)\}q^{t}$$
\noindent is a Hilbert modular form of weight 3/2.

\end{conj}

In \cite{YZZ}, the authors work by induction. To apply their method we would need to prove that the refined classes $\{Z(t,\varphi;H)\}$ are compatible with the tower of varieties attached to quadratic spaces  $V_x\hookrightarrow V$ of signature $(n,2)\times(n+1,1)^{r-1}\times(n+2,0)^{d-r}$ (in which case a generalization of \cite{KM} should imply that $\sum_{\substack{t\in\mathcal O_F\\ t\gg 0}}[Z(t,\varphi;H)]q^{t}$ is a Hilbert modular form of weight $\frac n2 +1$ with coefficients in $H^{r+1}(\mathrm{Sh}_H(G/Z,X)(\mathbf C),\mathbf C)$).

\subsection{A Gross-Kohnen-Zagier-type conjecture}

\paragraph{The Bruhat-Tits tree}
In this section we recall some basic facts about the Bruhat-Tits tree (see \cite{CJ} and \cite{Vi}). 

Let $v$ be a finite place of $F$. The vertices of the Bruhat-Tits tree of $\mathrm{PGL}_2(F_v)$ are the maximal orders of $\mathrm{M}_2(F_v)$. Such maximal orders are endomorphism rings of lattices in $F_v^2$ (\cite{Vi}, lemme 2.1). There is an oriented edge between two vertices $\mathcal O_1$ and $\mathcal O_2$ if and only if there exist $L_1,L_2$ lattices in $F_v^2$ such that $\mathcal O_i=\mathrm{End}(L_i)$, $L_2\subset L_1$ and $L_1/L_2\simeq \mathcal O_{F_v}/\varpi_v\mathcal O_{F_v}$. The intersection of the source and the target of paths of length $n$ correspond to level $v^n$ Eichler orders. 

Fix some quadratic extension $K/F$. This data allow us to organize the Bruhat-Tits tree. Let $\Psi : K_v\hookrightarrow \mathrm{M}_2(F_v)$ be a $F_v$-embedding of $K_v$. Let $\mathrm{M}_0(N)$ be the set of matrices in $\mathrm{M}_2(F_v)$ which are upper triangular modulo $N$. If $$\Psi(\mathcal O_{K_v})=\Psi(K_v)\cap \mathrm{M}_0(N),$$
\noindent we say that $\Psi$ has level $N$. We can organize the vertices of the tree in "levels", by privileging a direction. Each level corresponds to a level of embedding relatively to $\mathcal O_{K_v}$ i.e. to orders which are in the same orbit under $K_v^{\times}$. The maximal orders in $\mathrm{PGL}_2(F_v)$ which are maximally embedded are on the bottom of the tree.

Figures 2, 3 and 4 illustrate the dependence on the ramification type of $v$ in $K$.

\begin{figure}[!h]
\label{arbredecompose}
\centering
\includegraphics[width=13cm]{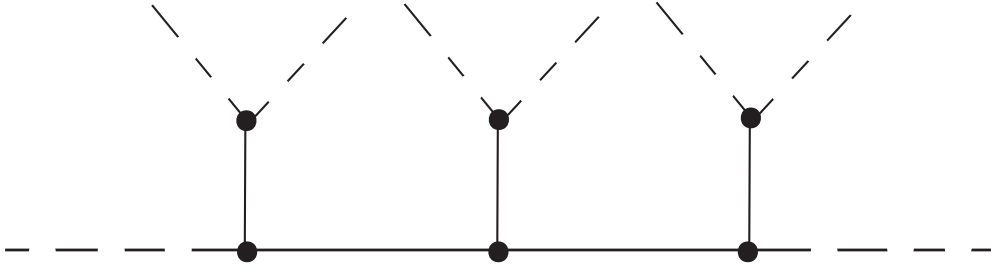}
\caption{Bruhat-Tits tree of $\mathrm{PGL}_2(F_v)$ when $v$ is split.} 
\end{figure}

\begin{figure}[!h]
\label{arbreramifie}
\centering
\includegraphics[width=13cm]{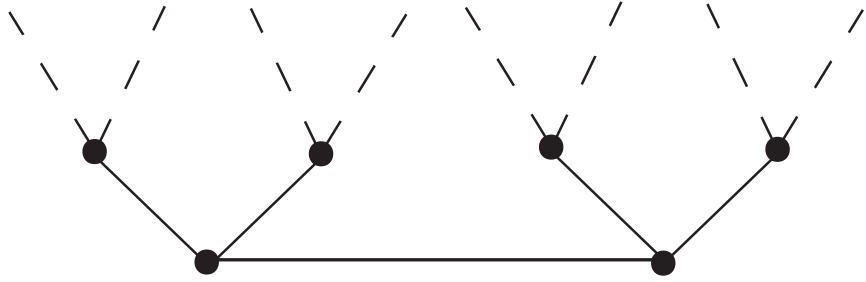}
\caption{Bruhat-Tits tree $\mathrm{PGL}_2(F_v)$ when $v$ is ramified.} 
\end{figure}

\begin{figure}[!h]
\label{arbreinerte}
\centering
\includegraphics[width=13cm]{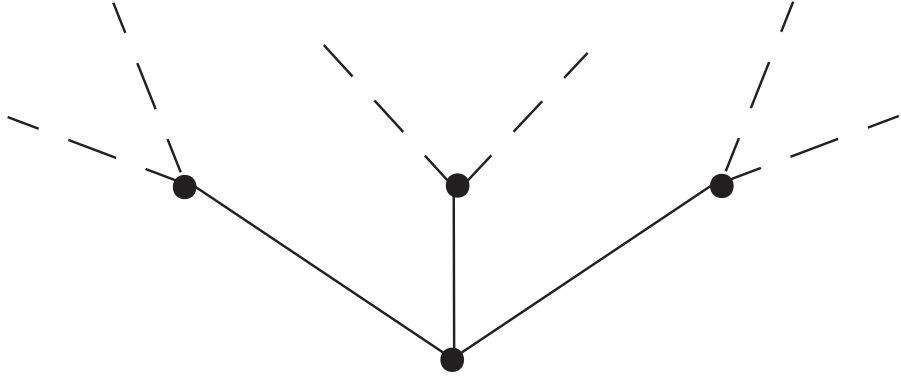}
\caption{Bruhat-Tits tree of $\mathrm{PGL}_2(F_v)$ when $v$ is inert.} 
\end{figure}

\paragraph{Darmon's points, Kudla's program and a Gross-Kohnen-Zagier-type theorem}

Recall that  $H=(R\otimes_{\mathbf Z}\widehat{\mathbf Z})^{\times}$, where $R$ is an Eichler order of $B$ of level $N_+$ and that $K=F+Fx$ satisfies the following Heegner hypothesis.

\begin{hyp}
\label{hypothesedeheegner}
 Each prime $\mathfrak p \mid N_+$ splits in $K$ and each prime $\mathfrak p \mid N_-$ is inert in $K$.
\end{hyp}

The group $G_x$ is isomorphic to $K^{\times}$ and $Z(x,1;H)$ is the image of $K^{\times}\backslash D_x\times \widehat{K}^{\times}/H$ in $\mathrm{Sh}_H(G,X)(\mathbf C)$. Note that

$$Z(x,1;H)=\mathscr T_1^1+t_1(\mathscr T_1^1),$$
\noindent where $\mathscr T_1^1=[\cup_{u\in\pi_0(T(\mathbf R))} q(u)\cdot \mathscr T^{\circ} ,1]_{H\widehat{F}^{\times}}$.

Let $\varphi=\mathbf{1}_{\widehat{R}^{\mathrm{Tr}=0}}.$ We are able to prove an analogue of Proposition A.I.1 of \cite{K1} when $N=1$, $B=\mathrm{M}_2(F)$, $R=\mathrm{M}_2(\mathcal O_F)$, $t=Q(x)=D_0\mathrm{nr}(x)\in F$ and $K=F+Fx$ is such that $K\cap R=\mathcal O_K$ and $\mathcal O_K=\mathcal O_F+\mathcal O_Fx$. Set $c_1(\mathscr T_1^1)=\{[t_1(x),b]_{H\widehat{F}^{\times}},\ b\in\widehat{B}^{\times}\}$.

\begin{prop}
\label{propkudlaappendice}
If $N=1$, $r=d$, $B=\mathrm{M}_2(F)$, $H=\widehat{R}^{\times}$ with $R=\mathrm{M}_2(\mathcal O_F)$ and if $\mathcal O_K=\mathcal O_F+\mathcal O_F x$, then $Z(t,\varphi;H)$ is equal to
$$Z(x,1;H)=\mathscr T_1^1+c_1(\mathscr T_1^1)=\mathscr T_1^1-\varepsilon\mathscr T_1^1.$$
\end{prop}

\begin{remark}

Under the strong hypotheses above, $\varepsilon=(-1)^d$ and the cycle obtained is zero when $d$ is even.

\end{remark}

\begin{proof}
By definition
$$Z(t,\varphi;H)=\sum_{g\in\widehat{K}^{\times}\backslash \widehat{B}^{\times}/\widehat{R}^{\times}}\mathbf{1}_{\widehat{R}^{\mathrm{Tr}=0}}(g^{-1}\cdot x) Z(x,g;H).$$
We have to determine $g\in\widehat{K}^{\times}\backslash \widehat{B}^{\times}/\widehat{R}^{\times}$ satisfying $g^{-1} x g\in \widehat{R}^{\mathrm{Tr}=0}$, i.e. $x\in g\widehat{R}^{\mathrm{Tr}=0}g^{-1}$. As $F^{\times}\subset K^{\times}$, 
$$ \widehat{K}^{\times}\backslash\widehat{B}^{\times}/\widehat{F}^{\times}\widehat{R}^{\times}=\prod_v{}'  K_v^{\times}\backslash B_v^{\times}/{R}_v^{\times}=\prod_v{}'  K_v^{\times}\backslash B_v^{\times}/F_v^{\times} {R}_v^{\times}.$$
\noindent This allows us to work locally with $K_v^{\times}\backslash B_v^{\times}/F_v^{\times} {R}_v^{\times}$, which is identified to the $K_v^{\times}$-orbits of maximal orders of $\mathrm{PGL}_2(F_v)$. This gives the following condition, $x_v\in g_v R_v g_v^{-1}.$

First let us consider those $g_v\in B_v^{\times}/R_v^{\times}F_v^{\times}$ satisfying $x_v\in g_vR_vg_v^{-1}$. The ring $g_vR_vg_v^{-1}$ is a maximal order containing $x_v$. Using the fact that $\mathcal O_K=\mathcal O_F+\mathcal O_F x$, we have
$$x_v\in g_vR_vg_v^{-1}\Longleftrightarrow g_vR_vg_v^{-1}\cap K_v=\mathcal O_{K_v}.$$
\noindent Hence the maximal order $g_vR_vg_v^{-1}$ is maximally embedded in $K_v$. It is identified to a vertex at the lowest level of the Bruhat-Tits tree. As each vertex at the same level is in the same  $K_v^{\times}$-orbit, we have
$$\forall v\qquad g_v=1\in K_v^{\times}\backslash B_v^{\times}/F_v^{\times}R_v^{\times}.$$
\noindent Thus $Z(t,\varphi;H)=Z(x,1;H)$
\noindent and as $D_{x_1}$ is a set of two points, $Z(x,1;H)$ is identified with $\mathscr T_1^1+c_1(\mathscr T_1^1)=\mathscr T_1^1-\varepsilon\mathscr T_1^1$, thanks to Proposition \ref{calculdessignes}.

\end{proof}

We now consider the case when $N=N_+N_-\neq 1$ is prime to $d_{K/F}$. The following proposition is true even if $B\neq \mathrm{M}_2(F)$ but we still assume that $R$ is an Eichler order of level $N_+$ and $\mathcal O_K=\mathcal O_F+\mathcal O_Fx$.

\begin{prop}
\label{propkudlaappendice2}
Let $N$ be the conductor of $E$. If $N$ is prime to $d_{K/F}$, then
 $$Z(t,\varphi;H)=\prod_{v\mid N}(1+\mathrm{inv}_v(B)\varepsilon_v)Z(x,1;H).$$
\end{prop}

\begin{proof}
The proof is analogous to the proof of Proposition \ref{propkudlaappendice}. Let us first compute the number of terms in $Z(t,\varphi;H)$. We need to determine for each $v$ the number of $K_v^{\times}$-orbits of oriented paths of length $\mathrm{ord}_v(N_+)$ in the Bruhat-Tits tree; this is equal to the number of $g_v$ such that $x_v\in g_vR_vg_v^{-1}.$

\begin{itemize}
\item If $v\nmid N$ then the same argument as in Proposition \ref{propkudlaappendice} shows that there is only one orbit.
\item If $v\mid N_-$,  $B_v$ is ramified and $v$ is inert in $K$. Hence $K_v^{\times}\backslash B_v^{\times}/R_v^{\times}F_v^{\times}=\{1,\pi_v\}$ where $\pi_v\in B_v^{\times}$ is an element whose reduced norm has order 1 at $v$; $\pi_v$ corresponds to the Atkin-Lehner involution. 
\item If $v\mid N_+$, $v$ splits in $K$. Denote by $v^\delta$ the level of the order $R_v$. Each Eichler order of level $v^\delta$ is the intersection of the origin and the target of an oriented path of length $\delta$. By hypothesis those orders are maximally embedded in $K_v$ and the path corresponding to $g_vR_vg_v^{-1}$ is contained in the lowest level of the tree. As $K_v^{\times}$ acts by translations on this level, there are exactly two $K_v^{\times}$-orbits corresponding to $g_v$ depending on the orientation. We have $g_v^+$ and $g_v^-$ which are exchanged by the Atkin-Lehner involution corresponding to $\begin{pmatrix} 0&\varpi_v\\ 1&0\end{pmatrix}$. 
\end{itemize}

Let $n$ be the number of prime ideals in the decomposition of $N$. The sum $Z(t,\varphi;H)$ has $2^n$ factors. Let $W$ be the sets of these factors. By definition $Z(x,g;H)=[\cdot g]Z(x,1;H)$. Using Proposition \ref{calculdessignes} we obtain
$$Z(t,\varphi;H)=\sum_{g\in W}[\cdot g]Z(x,1;H)=\prod_{v\mid N}(1+\mathrm{inv}_v(B)\varepsilon_v)Z(x,1;H).$$

\end{proof}

Let us conclude this paper by another conjecture. Assume 
 that $E(F)$ has rank 1. Denote by $P_0$ some generator of $E(F)$ modulo torsion. For each $t\in \mathcal O_F$ totally positive such that $(t)$ is square free and prime to $d_{K/F}$, denote by $K[t]$ the quadratic extension
$$K[t]=F(\sqrt{-D_0t}),$$
\noindent which satisfies the hypothesis used to build Darmon's points. Let $P_{t,1}$ be Darmon's point obtained for $K[t]$ and $b=1$, and set $$P_t=\mathrm{Tr}_{K[t]_1^+/F} P_{t,1}.$$
\noindent The point $P_t$ lies in $E(F)$ and there exists an integer $[P_t]\in\mathbf Z$ such that$$P_t=[P_t]P_0\ \mathrm{modulo\ torsion.}$$

Proposition \ref{propkudlaappendice2} together with Conjecture \ref{conjecturesixdeuxtrois} suggest the following (as in Conjecture 5.3 of \cite{DT}).

\begin{conj}
\label{lastconj}
There exists some Hilbert modular form $g$ of level 3/2 such that the $[P_t]$s are proportional to some Fourier coefficients of $g$.
\end{conj}

\begin{remark}

Using the analogy with the Gross-Kohnen-Zagier theorem, the integers $[P_t]$ should be (proportional to) square roots of $L(E_{-D_0t},1)$, where $E_{-D_0t}$ is the twist of $E$ by $-D_0t$.

\end{remark}

Let us end this paper with two open questions.

\begin{question}
Does Bruinier's genralization of Borcherds products \cite{Bruinier} give anything interesting in this situation ?

\end{question}

It is natural to expect that results of Cornut and Vatsal \cite{CV, CV2} hold also for Darmon's points.

\begin{question}

Would it be possible to deduce such a result from suitable equidistribution properties for the real tori $\mathscr T_b^{\circ}$ ?

\end{question}

\bibliography{bibliothese}{}
\bibliographystyle{alpha}

\end{document}